\documentclass[10pt,a4paper,fleqn]{article}

\usepackage[paper=a4paper,left=25mm,right=25mm,top=25mm,bottom=25mm]{geometry}
\usepackage[utf8]{inputenc}
\usepackage[english]{babel}
\usepackage[T1]{fontenc}
\usepackage{lmodern}
\usepackage{amsmath,amssymb,amsfonts,amsthm,mathtools}
\usepackage{mathrsfs}
\usepackage{bm}
\usepackage{booktabs}
\usepackage{enumitem}
\usepackage{graphicx}
\usepackage[numbers]{natbib}
\usepackage{microtype}
\usepackage{url}
\usepackage[hidelinks]{hyperref}
\usepackage{tikz}
\usepackage{pgfplots}
\pgfplotsset{compat=1.18}
\usepackage{comment}
\usepackage{longtable}

\allowdisplaybreaks

\newcommand{\R}{\mathbb{R}}
\newcommand{\N}{\mathbb{N}}
\newcommand{\E}{\mathbb{E}}
\newcommand{\Pp}{\mathbb{P}}
\newcommand{\cP}{\mathcal{P}}

\newcommand{\cE}{\mathcal{E}}
\newcommand{\cR}{\mathcal{R}}
\newcommand{\cW}{\mathcal{W}}
\newcommand{\Sph}{\mathbb{S}}
\newcommand{\Spd}{\mathbb{S}_{++}}

\newcommand{\BW}{\mathrm{BW}}
\newcommand{\tr}{\operatorname{tr}}
\newcommand{\Cov}{\operatorname{Cov}}
\newcommand{\law}{\operatorname{law}}

\newcommand{\Id}{\mathrm{Id}}

\newcommand{\dd}{\,\mathrm{d}}
\newcommand{\eqd}{\stackrel{\mathrm d}=}
\newcommand{\push}{_{\#}}
\newcommand{\norm}[1]{\left\lVert #1\right\rVert}
\newcommand{\abs}[1]{\left\lvert #1\right\rvert}

\newcommand{\trans}{\mathsf T}
\usepackage{xcolor}

\definecolor{light-gray}{gray}{0.95}
\definecolor{darkblue}{rgb}{0,0,.5}
\definecolor{foxred}{rgb}{0.7, 0.11, 0.11}

\newtheorem{theorem}{Theorem}[section]
\newtheorem{proposition}[theorem]{Proposition}
\newtheorem{corollary}[theorem]{Corollary}
\theoremstyle{definition}

\theoremstyle{plain}
\newtheorem{lemma}[theorem]{Lemma}
\theoremstyle{definition}

\theoremstyle{remark}
\newtheorem{remark}[theorem]{Remark}

\title{Quadratic Wasserstein transport for elliptical distributions
}

\author{
Jonathan Ansari\textsuperscript{1}
\qquad
Ludger R{\"u}schendorf\textsuperscript{2}\\[0.45em]
\small\textsuperscript{1} Paris Lodron Universität Salzburg, Austria\\
\small\texttt{jonathan.ansari@plus.ac.at}\\[0.35em]
\small\textsuperscript{2} Albert-Ludwigs-Universität Freiburg, Germany\\
\small\texttt{ruschen@stochastik.uni-freiburg.de}
}

\date{November 24, 2026}

\begin{document}
\maketitle

\begin{abstract}
The quadratic Wasserstein distance between Gaussian distributions is given by the classical Bures–Wasserstein formula, which depends only on the means, the covariance matrices \(\Sigma_0\) and \(\Sigma_1\), and their interaction \(K=(\Sigma_0^{1/2}\Sigma_1\Sigma_0^{1/2})^{1/2}\). 
The same formula holds for elliptical distributions with a common radial law. 
In this paper, we extend this representation to general elliptical distributions by identifying an additional transport cost \(\Gamma\) that depends on the radial laws and the spectrum of \(K\). 
We describe \(\Gamma\) via a weighted quadratic transport problem involving Dirichlet distributions and show that optimal couplings preserve directions within the eigenspaces of \(K\).
When \(K\) has a single eigenvalue, we obtain an explicit cost formula and an optimal Monge map.
For general spectra, we derive remarkably tight lower and upper bounds for \(\Gamma\) using nested spectral projections and Knothe--Rosenblatt couplings. The lower bound is based on one-dimensional integrals and can be computed fast. 
For Gaussian, Student-\(t\), uniform ball, and uniform sphere distributions, the upper bound also admits computationally efficient representations via dynamic programming.
Numerical experiments demonstrate that we reduce the gap between the previously known bounds by more than \(80\%\) in many settings, and by over \(99.9\%\) in certain scenarios, even in dimension \(10^6\). These improvements are particularly valuable in high-dimensional or heavy-tailed settings, where accurate numerical computation of quadratic Wasserstein distances is challenging.
\end{abstract}

\medskip



\section{Introduction}\label{sec:intro}

For probability measures $P_0,P_1\in\cP_2(\R^d)$, the quadratic Wasserstein distance is defined by
\begin{align}\label{eq:W2-def-intro}
W_2^2(P_0,P_1)
:=
\inf_{\pi\in\Pi(P_0,P_1)}
\int_{\R^d\times\R^d}\norm{x-y}^2\,\pi(\dd x,\dd y),
\end{align}
where $\Pi(P_0,P_1)$ denotes the set of couplings of $P_0$ and $P_1$. For Gaussian distributions
$P_i=N(m_i,\Sigma_i)$, $i=0,1$, one has the classical formula
\begin{equation}\label{eq:gaussian-formula-intro}
W_2^2(P_0,P_1)
=
\norm{m_0-m_1}^2+d_{\BW}^2(\Sigma_0,\Sigma_1),
\end{equation}
with the \emph{Bures-Wasserstein distance}
\begin{align}\label{eq:bures-def-intro}
d_{\BW}^2(\Sigma_0,\Sigma_1)
:=
\tr(\Sigma_0)+\tr(\Sigma_1)
-2\tr\!\left[\left(\Sigma_0^{1/2}\Sigma_1\Sigma_0^{1/2}\right)^{1/2}\right].
\end{align}
The classical formula for Gaussian distributions was obtained in several equivalent forms in the early optimal transport and matrix literature; see, among others, \cite{DowsonLandau1982,OlkinPukelsheim1982,GivensShortt1984,Ruschendorf1985}. Its geometry was subsequently developed by Takatsu \cite{Takatsu2011}. For a matrix-analytic account of the Bures--Wasserstein metric, see \cite{BhatiaJainLim2019}; 
for a recent extension to adapted transport between Gaussian processes, we refer to \citet{GunasingamWong2025}.

Two fundamental extensions of the Gaussian formula are due to \citet{Gelbrich1990}. First, for arbitrary probability distributions
\(P_0,P_1\in\cP_2(\R^d)\) with means \(m_0,m_1\) and covariance matrices
\(\Sigma_0,\Sigma_1\), the \emph{Gelbrich bound} states that
\begin{equation}\label{eq:gelbrich-bound-intro}
W_2^2(P_0,P_1)
\geq
\norm{m_0-m_1}^2
+
d_{\BW}^2(\Sigma_0,\Sigma_1).
\end{equation}
Second, \citet{Gelbrich1990} showed that equality holds for elliptical distributions
with a common standardized radial law, i.e.,
\begin{align}\label{eq_Gelb}
    W_2^2(P_0,P_1)
=
\norm{m_0-m_1}^2
+
d_{\BW}^2(\Sigma_0,\Sigma_1) \qquad \text{for } P_i = \cE_2(\mu_i,\Sigma_i,\nu),
\end{align}
where $\cE_2(\mu,\Sigma,\nu)$ denotes the square-integrable
elliptical distribution with mean $\mu\in\R^d$,
positive definite covariance matrix
$\Sigma\in\R^{d\times d}$, and standardized radial
distribution $\nu$ on $[0,\infty)$ satisfying
\(\int_0^\infty r^2\,\nu(\mathrm dr)=d\).
Recall that any such distribution admits a stochastic
representation
\begin{align}\label{eq_rep_ell}
X=\mu+RAU,
\qquad
R\sim\nu,\quad
U\sim\mathcal U(\Sph^{d-1}),\quad
R\perp U,\quad
AA^\top=\Sigma,
\end{align}
where \(\mathcal{U}(\Sph^{d-1})\) denotes the uniform distribution on the unit sphere of \(R^d\).
Note that, since $\mathbb E[UU^\top]=I_d/d$, the
normalization ensures $\operatorname{Cov}(X)=\Sigma$.


The situation changes substantially when the radial laws are different. Radial transports between spherically equivalent distributions are classical; see \citet{CuestaRuschendorfTuero1993}. More recently, such transports have been used for Wasserstein barycenters by \citet{GhaffariWalker2023}. Further, \citet{ChenZhang2024} study radial-contoured distributions and emphasize the lack of explicit formulas for general elliptical distributions with different generators. \citet{GhaffariWalker2021} obtain explicit maps for classes related by common copulas and orthonormal transformations, including a class of scale mixtures of normals. Related lower and upper bounds for affine transformations of random vectors are developed by \citet{HammKorzeniowski2024}. A key difficulty is that composing an optimal radial transport with an optimal positive definite linear transport does not generally yield an optimal coupling, as shown in \cite[Proposition~3.14]{CuestaRuschendorfTuero1993}.

The challenge of determining the Wasserstein distance between elliptical distributions with different radial laws was also highlighted by Wong in the collection of open problems in \cite[Section~8, Question~2]{SeiMatsuzoe2025}.
More precisely, Wong
points out the absence of a general closed formula for different radial laws, provides a radial upper bound, and asks for an improvement of the Gelbrich lower bound that also depends on the radial laws. 

\subsubsection*{Spectral reduction}

The aim of the present paper is to determine the additional
transport cost arising from different radial laws of
elliptical distributions. As the first main result of this paper, we provide an 
exact decomposition of this cost based on the spectral
structure of the covariance matrices. To this end, we consider 
for $P_i=\cE_2(m_i,\Sigma_i,\nu_i)$, $i=0,1$,
the positive definite \emph{Bures matrix}
\begin{align}\label{eq:K-main-intro}
K:=
\left(
\Sigma_0^{1/2}\Sigma_1\Sigma_0^{1/2}
\right)^{1/2}.
\end{align}
This matrix $K$ already occurs in the Bures-Wasserstein distance \eqref{eq:bures-def-intro} and arises from the Bures-aligned covariance
factors \(A_0=\Sigma_0^{1/2}\) and \(A_1=\Sigma_0^{-1/2}K\),
which satisfy
$A_iA_i^{\trans}=\Sigma_i$ and
$A_0^{\trans}A_1=K$.
Thus, the transformations
$A_i^{-1}(X_i-m_i)$ yield standardized spherical
distributions \(\eta_{\nu_i} := \law(R_iU_i)\), while $K$ captures the relative
covariance geometry of the two elliptical laws; see also Figure \ref{fig:optimal-monge-transport}.

Now, let $\bm\kappa=(\kappa_1,\ldots,\kappa_k)$ denote
the distinct eigenvalues of $K$, with corresponding
eigenspace dimensions $n_1,\ldots,n_k$.
Our first main result, Theorem \ref{thm:main}, establishes the representation
\begin{align}\label{eq:main-result-intro}
W_2^2(P_0,P_1)
=
\norm{m_0-m_1}^2
+d_{\BW}^2(\Sigma_0,\Sigma_1)
+\cW_{2,\bm\kappa}^2(\beta_0,\beta_1),
\end{align}
where $\cW_{2,\bm\kappa}$ denotes the weighted quadratic Wasserstein distance in \eqref{eq:spherical-spectral-reduction-intro} between the \(k\)-dimensional distributions \(\beta_0\) and \(\beta_1\). The latter are the joint distributions of the block radii, which are the radial parts within the pairwise orthogonal eigenspaces, i.e.,
\begin{align*}
\beta_i = \law \big(
R_i\sqrt{D_1},\ldots,R_i\sqrt{D_k}\big), \quad R_i\sim\nu_i \text{ independent of } (D_1,\ldots,D_k)
\sim
\mathrm{Dirichlet}
\left(\tfrac{n_1}{2},\ldots,\tfrac{n_k}{2}\right);
\end{align*}
see Proposition \ref{prop:spherical-reduction-intro}.
Consequently, the original $d$-dimensional
transport problem is reduced to a weighted
transport problem whose dimension is the number \(k\)
of distinct eigenvalues of $K$.

\subsubsection*{Generator effect}

The last term in \eqref{eq:main-result-intro}
describes the additional transport cost
caused by different standardized radial laws.
We call this quantity the \emph{generator effect}, defined by
\begin{equation}\label{eq:generator-effect-intro}
\begin{aligned}
\Gamma(P_0,P_1)
&:= \cW_{2,\bm\kappa}^2(\beta_0,\beta_1) = W_2^2(P_0,P_1)
-\norm{m_0-m_1}^2
-d_{\BW}^2(\Sigma_0,\Sigma_1).
\end{aligned}
\end{equation}
Since $K$ is positive definite, the generator
effect is positive if and only if
$\nu_0\neq\nu_1$. It is important to note that this quantity depends not only on the radial laws but also on the covariance geometry encoded by $K$.

\begin{figure}[t!]
    \centering
    \includegraphics[
        width=0.60\textwidth,
        keepaspectratio
    ]{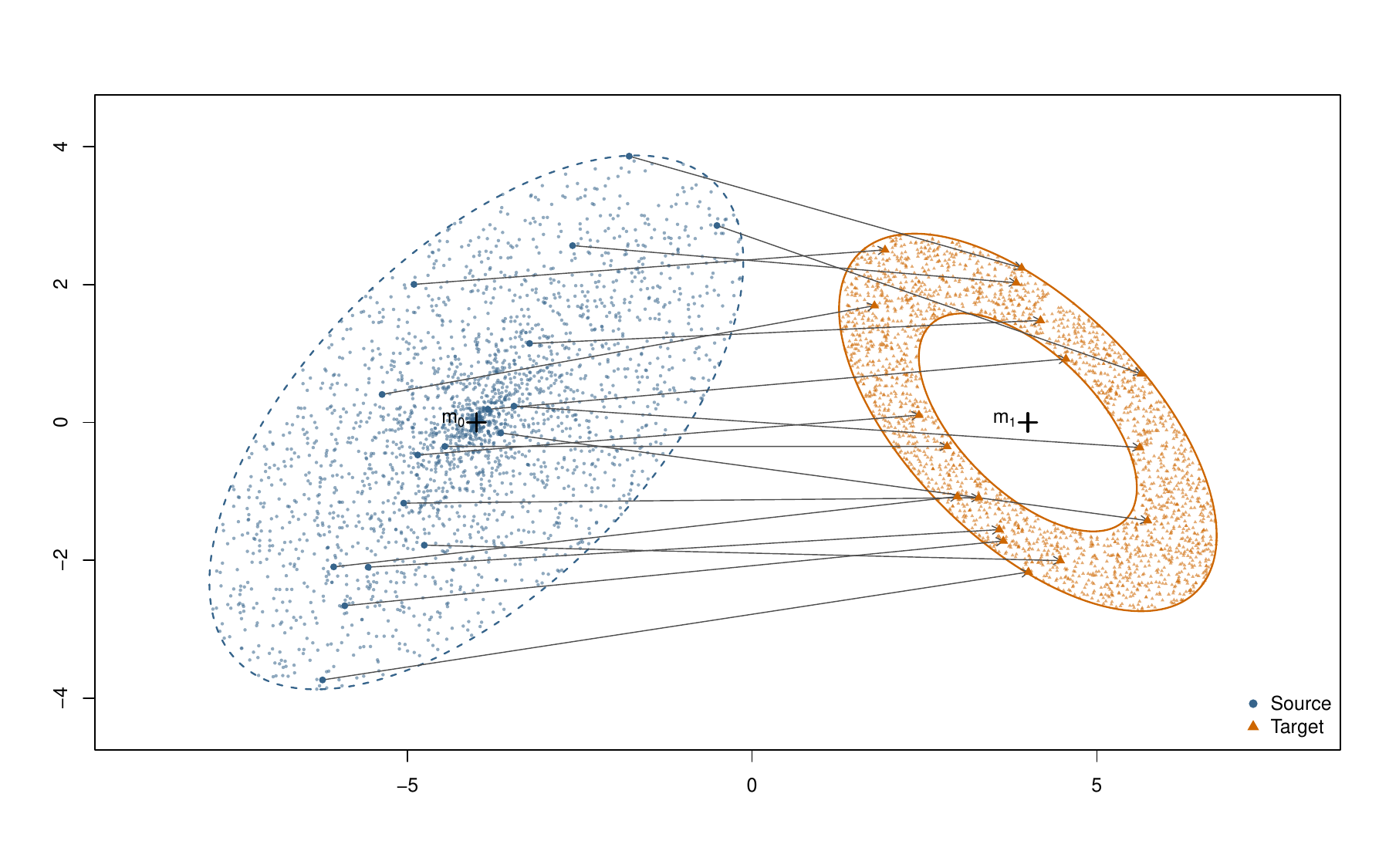}
    \hfill
    \includegraphics[
        width=0.38\textwidth,
        keepaspectratio
    ]{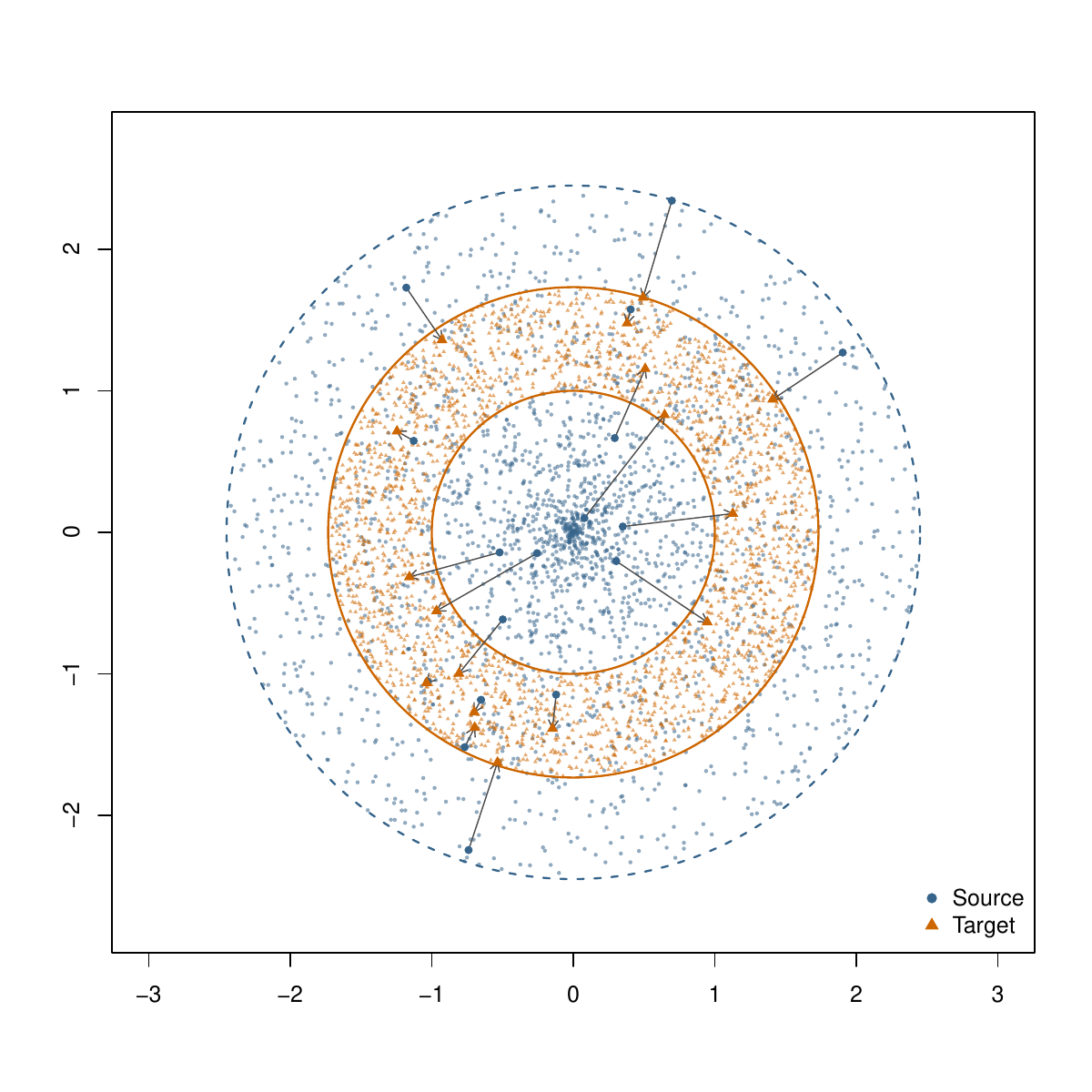}

\caption{Optimal Monge transport between the elliptical distributions
$P_i=\cE_2(m_i,\Sigma_i,\nu_i)$, $i=0,1$,
with means $m_0=(-4,0)^\top$ and $m_1=(4,0)^\top$,
covariance factors
$A=\left(\begin{smallmatrix}3/2&1/2\\1/2&3/2\end{smallmatrix}\right)$
and
$B=\left(\begin{smallmatrix}3/2&-1/2\\-1/2&3/2\end{smallmatrix}\right)$,
and radial variables 
$R_0\sim \nu_0 =\mathcal{U}(0,\sqrt{6})$ and
$R_1\sim \nu_1 = \sqrt{\mathcal{U}(1,3)}$.
The covariance matrices satisfy
$\Sigma_0=AA^\top$ and $\Sigma_1=BB^\top$,
while the Bures matrix is
$K=A^\top B=2I_2$. 
Left: Optimal transport in the original elliptical coordinates.
Right: The corresponding transport in standardized spherical
coordinates, where the optimal radial map is given by
$\tau(r)= F_{\nu_1}^{-1}(F_{\nu_0}(r)) = \sqrt{1+2r/\sqrt{6}}$.
Arrows connect points matched by the optimal Monge map
$T$ in Corollary~\ref{cor:scalar-intro} with cost \(W_2^2(P_0,P_1) =
64+2+\frac{120-72\sqrt2-4\sqrt6}{15} \approx 66.56\).
}
\label{fig:optimal-monge-transport}
\end{figure}

In the \emph{isotropic} case, where \(K = \lambda I_d\) has only one eigenvalue \(\lambda\) (i.e., \(k=1\)), we can determine the transport cost \(W_2(P_0,P_1)\) and the optimal Monge map from \(P_0\) to \(P_1\) explicitly; see Corollary \ref{cor:scalar-intro} and Figure \ref{fig:optimal-monge-transport}.
However, in the \emph{anisotropic} case (i.e., \(k\geq 2\)), there is generally no analogous closed-form expression available since the spectral reduction leads to the multidimensional optimal transport problem \(\cW_{2,\bm\kappa}(\beta_0,\beta_1)\) involving the joint distributions of the block radii.

\subsubsection*{Improved bounds}

Since no explicit formula for \(\cW_{2,\bm\kappa}(\beta_0,\beta_1)\) is available, computable lower and upper bounds are of particular interest. This question was also raised by Wong in
\cite[Section~8, Question~2]{SeiMatsuzoe2025},
who provides the elementary upper bound
\begin{align}\label{eq:wong-bound-intro}
U_{\mathrm{Wong}}
:=
\frac{\tr K}{d}
W_2^2(\nu_0,\nu_1),
\end{align}
which follows by coupling the radial variables
comonotonically and choosing a common uniform
spherical direction, independent of the radii.
Wong also asks for a radial-sensitive improvement of Gelbrich's lower bound \(L_{\mathrm{Gelbrich}} := 0\) obtained from \eqref{eq:gelbrich-bound-intro} and \eqref{eq_Gelb}.
In Corollary~\ref{cor:intro-radial-bounds}, we derive
an elementary lower bound \(\lambda_{\min}(K)\delta_R^2\) from the smallest eigenvalue
of $K$ and from the optimal radial transport between the
standardized spherical distributions. In general, however, these elementary bounds are too wide to provide an accurate estimate of the generator effect.

As a second main contribution of our paper, we substantially
tighten the elementary lower and upper bounds for the generator effect
$\Gamma(P_0,P_1)$ through
\begin{align}\label{eq_bound_geneff}
    0= L_{\mathrm{Gelbrich}} \leq \lambda_{\min}(K) \delta_R^2 \leq L_K \leq \Gamma(P_0,P_1) \leq U_K \leq U_{\mathrm{Wong}}.
\end{align}
Our lower bound $L_K$ in Theorem~\ref{thm:nested-bound}
is based on nested spectral projections. More precisely,
we consider projections onto nested sums of eigenspaces
of $K$ corresponding to increasingly large eigenvalues.
This construction exploits the spectral structure of $K$
and yields an explicit lower bound in terms of
one-dimensional Wasserstein distances between projected
radial laws.

Our upper bound $U_K$ in
Theorem~\ref{thm:multi-conditional-upper}
is obtained by coupling block-radius vectors through
conditionally comonotone Knothe--Rosenblatt couplings.
In general, evaluating the corresponding transport costs
requires multidimensional integration. However, as shown
in Theorem~\ref{lem:KR-student-one-dimensional},
for arbitrary pairs of Gaussian, Student-$t$, uniform
ball, and uniform sphere distributions, these costs
can be expressed in terms of one-dimensional integrals.
Moreover, the construction of \(U_K\) is based on an optimization over all combinations of projections into orthogonal eigenspaces of the matrix \(K\). In Theorem~\ref{lem:KR-student-one-dimensional}, we speed up this optimization via dynamic programming and show that it can be carried out efficiently with $O(d^2)$ transitions.

Our numerical experiments demonstrate that the resulting
bounds are remarkably tight, even in dimensions as
large as $d=10^6$. 
Compared with the classical Gelbrich--Wong bounds, our bounds
reduce the width of the corresponding interval by more than
$80\%$ in many settings and by over $99.9\%$ in some cases;
see Table~\ref{tab:intro-highlights} and
Section~\ref{sec:simulation} for details.
Such tight bounds are particularly valuable in high-dimensional
or heavy-tailed settings, where accurate numerical computation
of quadratic optimal transport costs is challenging.
Sample-based methods suffer from the curse of dimensionality,
with convergence rates also depending on the tail behavior
of the distributions \cite{FournierGuillin2015,WeedBach2019},
while entropically regularized methods such as Sinkhorn
introduce an additional approximation bias
\cite{Cuturi2013,GenevayEtAl2019}.
By exploiting the spectral structure of elliptical
distributions, our approach provides direct lower and upper
bounds for the unregularized Wasserstein transport cost
without requiring high-dimensional discretization or
entropic regularization.

Such computationally tractable bounds are relevant to
applications of Wasserstein geometry, including elliptical
distributional embeddings \cite{MuzellecCuturi2018},
statistical models of affine deformations
\cite{AmariMatsuda2024}, and Wasserstein barycenters
and mixtures of radially contoured distributions
\cite{GhaffariWalker2023,ChenZhang2024,ChenWangZhang2026}.
In these settings, our results facilitate the comparison
of models with different tail behavior.
They may also be useful for sensitivity analysis and
the calibration of Wasserstein ambiguity sets in
distributionally robust optimization \cite{KuhnEtAl2019}.

\begin{table}[t!]
\centering\scriptsize
\setlength{\tabcolsep}{2.1pt}
\renewcommand{\arraystretch}{1.08}
\begin{tabular}{@{}rllrrrrrrr@{}}
\toprule
$d$ & Radial pair & \shortstack{Bures\\matrix} & $k$ & $L_{\min}$ & $L_K$ & $\Gamma(P_0,P_1)$ & $U_K$ & $U_{\mathrm{Wong}}$ & Improvement \\
\midrule
$2$ & $N$--$t_{3}$ & $K_{\mathrm{ar}}$ & $2$ & $0.5325$ & $0.773$ & $\approx 0.79$ & $0.7987$ & $0.7987$ & $96.78\%$ \\
$2$ & $N$--$t_{3}$ & $K_{\mathrm{pow}}$ & $2$ & $0.3195$ & $0.7524$ & $\approx 0.77$ & $0.784$ & $0.7987$ & $96.04\%$ \\
$2$ & $N$--$t_{3}$ & $K_2$ & $2$ & $0.3994$ & $0.7601$ & $\approx 0.78$ & $0.7987$ & $0.7987$ & $95.17\%$ \\
\addlinespace[7pt]
$10$ & $N$--$t_{3}$ & $K_2$ & $2$ & $13.2747$ & $23.2587$ & $\approx 24$ & $24.8822$ & $26.5494$ & $93.89\%$ \\
$10$ & $N$--$t_{3}$ & $K_{\mathrm{ar}}$ & $10$ & $4.8272$ & $24.6431$ & -- & $26.1541$ & $26.5494$ & $94.31\%$ \\
$10$ & $N$--$t_{3}$ & $K_{\mathrm{pow}}$ & $10$ & $1.7131$ & $20.2948$ & -- & $21.9639$ & $26.5494$ & $93.71\%$ \\
\addlinespace[7pt]
$100$ & $N$--$t_{3}$ & $K_2$ & $2$ & $1\,738.367$ & $3\,212.9691$ & $\approx 3\,300$ & $3\,339.4969$ & $3\,476.7339$ & $96.36\%$ \\
$100$ & $N$--$t_{3}$ & $K_{\mathrm{ar}}$ & $100$ & $68.8462$ & $3\,329.8786$ & -- & $3\,437.772$ & $3\,476.7339$ & $96.90\%$ \\
$100$ & $N$--$t_{3}$ & $K_{\mathrm{pow}}$ & $100$ & $21.2646$ & $2\,098.3883$ & -- & $2\,280.9218$ & $3\,476.7339$ & $94.75\%$ \\
\addlinespace[7pt]
$10^6$ & $N$--$t_{3}$ & $K_2$ & $2$ & $2.0180\times10^{11}$ & $4.032050\times10^{11}$ & $\approx 4.0334\times10^{11}$ & ${4.033801\times10^{11}}^{\dagger}$ & $4.0360\times10^{11}$ & $99.956610\%$ \\
$10^6$ & $N$--$t_{3}$ & $K_{\mathrm{ar}}$ & $10^6$ & $807\,191.6524$ & ${4.033847\times10^{11}}^{\dagger}$ & -- & ${4.035401\times10^{11}}^{\dagger}$ & $4.0360\times10^{11}$ & $99.961480\%$ \\
$10^6$ & $N$--$t_{3}$ & $K_{\mathrm{pow}}$ & $10^6$ & $245\,357.2355$ & ${2.1089\times10^{11}}^{\dagger}$ & -- & ${2.2914\times10^{11}}^{\dagger}$ & $4.0360\times10^{11}$ & $95.48\%$ \\
\bottomrule
\end{tabular}
\caption{Bounds in \eqref{eq_bound_geneff} for the generator effect $\Gamma(P_0,P_1)$ between Gaussian and Student-$t_3$ laws. The spectra $K_2,K_{\mathrm{ar}},K_{\mathrm{pow}}$ are specified in Table~\ref{tab:spectra-simple} and satisfy $\operatorname{tr}K=d^2$. A dagger marks the value for $L^-$ or $U^+$ at $d=10^6$; see \eqref{eq:simulation-conservative-bounds}. A dash denotes that no numerical approximation is reported for $\Gamma(P_0,P_1)$; see Section~\ref{sec:simulation} for details.}
\label{tab:intro-highlights}
\end{table}

\subsubsection*{Organization of the paper}

The remainder of the paper is organized as follows.
Section~\ref{sec_22} introduces the mathematical framework,
including the Bures decomposition and the spectral
reduction to block-radius distributions.
Section~\ref{sec_main} presents our main results: the exact
representation of the quadratic Wasserstein distance,
the explicit optimal Monge map in the isotropic case,
and the improved lower and upper bounds for the
generator effect.
Section~\ref{sec:dyn} develops an efficient computation of the
upper bound using one-dimensional integrals and
dynamic programming.
Section~\ref{sec:simulation} provides numerical experiments illustrating
the accuracy and computational feasibility of our
bounds, including high-dimensional settings.
Some auxiliary results and all proofs are collected
in the Appendix.

\section{Mathematical framework}\label{sec_22}

Throughout, $d\geq1$ is fixed and $\R_+=[0,\infty)$. Further, $\cP_2(\R^d)$ denotes the Borel probability measures on $\R^d$ with finite second moment. For $P,Q\in\cP_2(\R^d)$, $\Pi(P,Q)$ is the set of Borel probability measures on $\R^d\times\R^d$ with marginals $P$ and $Q$. The quadratic Wasserstein distance is defined in \eqref{eq:W2-def-intro}. We use standard existence results for optimal couplings without further comment; see, for example, Rachev and R\"uschendorf \cite{RachevRuschendorf1998} and Villani \cite{Villani2009}. 

In the context of elliptical distributions with second moments, it is convenient to consider a normalized parametrization. Therefore, let
\begin{equation}\label{eq:Rclass-intro}
\cR_d
:=
\left\{
\nu\in\cP_2(\R_+):
\int_0^\infty r^2\,\nu(\mathrm{d} r)=d
\right\}
\end{equation}
denote the class of standardized radial laws and \(\Spd^d\) the set of positive definite matrices with values in \(\R^{d\times d}\). For \(m\in\R^d\),
\(\Sigma\in\Spd^d\), and \emph{generator} \(\nu\in\cR_d\), recall 
\begin{equation}\label{eq:elliptical-def-intro}
\cE_2(m,\Sigma,\nu)
:=
\law\!\left(m+ARU\right),
\qquad
AA^{\trans}=\Sigma,
\end{equation}
where \(R\sim\nu\), \(U\) is uniform on the Euclidean unit sphere
\(\Sph^{d-1}\), and \(R\) and \(U\) are independent. The law in
\eqref{eq:elliptical-def-intro} does not depend on the particular factor
\(A\). Moreover, the square-integrability and normalization in \eqref{eq:Rclass-intro} ensure that \(m\) and \(\Sigma\) are the mean and covariance matrix; see Lemma \ref{lem:ell-normalization}.
For a measurable map $T$, the push-forward $T\push P$ is defined by $(T\push P)(B)=P(T^{-1}(B))$. The notation $\law(X)$ denotes the distribution of $X$ and $X\eqd Y$ equality in distribution. For a finite interval \(J\subset \R\), we write \(\mathcal{U}(J)\) for the uniform distribution on \(J\).
For a subspace \(E\subseteq\R^d\), \(I_E\) denotes the identity operator on \(E\),
while \(\Id\) denotes the identity map on the ambient space under consideration.

In the remainder of this section, we introduce
the concepts and preliminary results needed
for our main results in Section~\ref{sec_main}.

\subsection{Preliminaries on quadratic optimal transport}

First, we recall a fundamental theorem in quadratic optimal
transport that characterizes optimal Monge maps as the subgradient of a convex function.  

\begin{theorem}[Quadratic optimality and convex subgradients, \cite{RuschendorfRachev1990,Brenier1991}]\label{lem:quadratic-optimality}
Let $P,Q\in\cP_2(\R^d)$. A coupling $\pi\in\Pi(P,Q)$ is optimal with respect to the quadratic Wasserstein distance if and only if there is a proper lower semicontinuous convex function $\phi$ such that
\begin{equation}\label{eq:subgradient-optimality}
 y\in\partial\phi(x)\quad\text{for }\pi\text{-almost every }(x,y),
\end{equation}
where $\partial\phi(x)$ denotes the convex subdifferential of $\phi$ at $x$.  In particular, if a measurable map $T$ satisfies $T\push P=Q$ and $T(x)\in\partial\phi(x)$ for $P$-almost every $x$, then $(\Id,T)\push P$ is an optimal Monge coupling. If $P$ is absolutely continuous, the optimal coupling is unique and is induced by a convex-gradient map $T=\nabla\phi$.
\end{theorem}

The existence and subgradient characterization is due to R\"uschendorf and
Rachev \cite{RuschendorfRachev1990}, obtained by duality and convex analysis. In the Lebesgue-continuous case this characterization implies in particular the existence of an optimal convex gradient transport map.
For existence and uniqueness of convex-gradient transport maps under
Lebesgue-absolute continuity of the source, see
Brenier \cite{Brenier1987,Brenier1991}.

The following classical identity for spherical laws is an application of the convex gradient characterization in Theorem \ref{lem:quadratic-optimality}. It is obtained by coupling the radii comonotonically and using a common uniform spherical direction independent of the radii; for an illustration, see the right-hand plot in Figure \ref{fig:optimal-monge-transport}. Recall that, for
\(\nu\in\cR_d\), the associated standardized spherical law is denoted by
\begin{align}\label{def:eta_nu}
    \eta_\nu:=\law(RU), \qquad R\sim\nu, ~ U\sim\mathcal{U}(\Sph^{d-1}), ~ R\perp U.
\end{align}
Note also that spherical laws are probability distributions on \(\R^d\) that are invariant under orthogonal transformations.

\begin{lemma}[Optimal transport between spherical laws, {\cite{CuestaRuschendorfTuero1993}}]\label{lem:radial-optimal}
Let $\nu_0,\nu_1\in\cR_d$.  Then we have
\begin{align}
    W_2^2(\eta_{\nu_0},\eta_{\nu_1})
=
W_2^2(\nu_0,\nu_1)
= 
\int_0^1
\left|
F_{\nu_0}^{-1}(u)-F_{\nu_1}^{-1}(u)
\right|^2
\,\dd u.
\end{align}
If $\nu_0$ is non-atomic, the optimal Monge map for \(W_2(\eta_{\nu_0},\eta_{\nu_1})\) is
\begin{equation}\label{eq:radial-map-standard}
S(z)
=
\begin{cases}
\displaystyle\frac{\tau(\norm z)}{\norm z}z,&z\neq0,\\[0.6em]
0,&z=0,
\end{cases}
\end{equation}
where $\tau=F_{\nu_1}^{-1}\circ F_{\nu_0}$ pushes $\nu_0$ to $\nu_1$. 
\end{lemma}

For general elliptical distributions, which are just linear transformations of spherical distributions according to \eqref{eq_rep_ell}, an optimal coupling exists. If the source distribution is
absolutely continuous, Theorem~\ref{lem:quadratic-optimality} further guarantees a unique optimal Monge map. Its
explicit form is generally unknown.

In the following, we exploit the spectral decomposition of $K$ to reduce the dimension of the underlying optimal transport problem. When $K$ has only one eigenvalue, we still obtain an explicit formula for the transport cost and an optimal Monge
map. More generally, we use the spectral structure of $K$ to derive
remarkably sharp lower and upper bounds for the optimal transport cost in various settings.

\subsection{Bures alignment and exact cost decomposition}

Now, let \(\nu_i\in\cR_d\) be arbitrary and consider \(
P_i=\cE_2(m_i,\Sigma_i,\nu_i)\)
with \(m_i\in\R^d\) and \(\Sigma_i\in\Spd^d\), \(i=0,1\).
The decomposition of \(\cE_2(m_i,\Sigma_i,\nu_i)\) into location, covariance geometry, and radial part is not merely a parametrization of elliptical distributions. For the quadratic cost, it is also closely related to an exact
couplingwise identity. To see this, let $X_0\sim P_0$ and $X_1\sim P_1$. Since
$\E(X_0-m_0)=\E(X_1-m_1)=0$, every coupling satisfies
\begin{equation}\label{eq:mean-separation-intro}
 \E\norm{X_0-X_1}^2
 =\norm{m_0-m_1}^2
  +\E\norm{(X_0-m_0)-(X_1-m_1)}^2.
\end{equation}
Thus the location term separates before any optimization is performed.

We next separate the covariance geometry. Therefore, recall that the positive definite Bures matrix \(K\) in \eqref{eq:K-main-intro} satisfies 
$A_0=\Sigma_0^{1/2}$ and $A_1:=\Sigma_0^{-1/2}K$ as well as $A_0A_0^{\trans}=\Sigma_0$, $A_1A_1^{\trans}=\Sigma_1$, and $A_0^{\trans}A_1=K$. Hence, after centering, we may write
\begin{equation}\label{eq:aligned-intro}
 X_0-m_0=A_0Z_0,
 \qquad
 X_1-m_1=A_1Z_1,
\end{equation}
where $Z_i=R_iU_i\sim\eta_{\nu_i}$ are standardized
spherical random vectors. Recall from \eqref{def:eta_nu} that \(\eta_{\nu_i}=\law(R_iU_i)\) is centered and has covariance
$I_d$.

The choice of \(A_0\) and \(A_1\) makes the covariance contribution separate from the remaining transport cost. Indeed, for every coupling of \(Z_0\) and \(Z_1\), the following two quadratic expansions have the same coupling-dependent cross term:
\begin{align*}
 \E\norm{A_0Z_0-A_1Z_1}^2
 &=\tr(\Sigma_0)+\tr(\Sigma_1)-2\E[Z_0^{\trans}KZ_1], \\
 \E[(Z_0-Z_1)^{\trans}K(Z_0-Z_1)]
 &=2\tr(K)-2\E[Z_0^{\trans}K Z_1].
\end{align*}
Comparing these two identities and using the definition of the Bures--Wasserstein distance in \eqref{eq:bures-def-intro} gives
\begin{equation}\label{eq:couplingwise-preview-intro}
 \E\norm{A_0Z_0-A_1Z_1}^2
 =d_{\BW}^2(\Sigma_0,\Sigma_1)
  +\E\bigl[(Z_0-Z_1)^{\trans}K(Z_0-Z_1)\bigr].
\end{equation}
This identity holds \emph{before} taking an infimum. Together with
\eqref{eq:mean-separation-intro}, it shows that the mean contribution and the
Bures covariance contribution separate, while all remaining freedom
is contained in a weighted transport problem between the standardized spherical laws.

For \(M\in\Spd^d\) and \(P,Q\in\cP_2(\R^d)\), we define the weighted quadratic
Wasserstein distance by
\begin{equation}\label{eq:W2M-intro}
W_{2,M}^2(P,Q)
:=
\inf_{\gamma\in\Pi(P,Q)}
\int_{\R^d\times\R^d}
(z-w)^{\trans}M(z-w)\,\gamma(\dd z,\dd w).
\end{equation}
This gives
\begin{align*}
    W_{2,M}^2(P,Q)
=
W_2^2\bigl((M^{1/2})\push P,(M^{1/2})\push Q\bigr),
\end{align*}
and, since \(M\) is positive definite, \(W_{2,M}\) is a metric on
\(\cP_2(\R^d)\). We summarize the preceding decomposition of the transport cost as follows.

\begin{proposition}[Transport cost decomposition]\label{lem:effect-intro}
For \(P_i=\cE_2(m_i,\Sigma_i,\nu_i)\), \(i=0,1\), let \(K\) be given by
\eqref{eq:K-main-intro}. Then
\begin{equation}\label{eq:effect-formula}
W_2^2(P_0,P_1)
=
\norm{m_0-m_1}^2
+d_{\BW}^2(\Sigma_0,\Sigma_1)
+W_{2,K}^2(\eta_{\nu_0},\eta_{\nu_1}),
\end{equation}
where
\begin{equation}\label{eq:effect-equivalent}
W_{2,K}^2(\eta_{\nu_0},\eta_{\nu_1})
=
W_2^2\!\left(
\cE_2(0,K,\nu_0),
\cE_2(0,K,\nu_1)
\right).
\end{equation}
\end{proposition}

\begin{remark}[Generator effect]
\label{rem:effect-interpretation}
\begin{enumerate}[label=(\alph*)]

\item Gelbrich's inequality \eqref{eq:gelbrich-bound-intro} states that the quantity
\begin{align*}
  W_{2,K}^2(\eta_{\nu_0},\eta_{\nu_1}) =  W_2^2(P_0,P_1)
-
\norm{m_0-m_1}^2
- d_{\BW}^2(\Sigma_0,\Sigma_1)
\end{align*}
is non-negative and vanishes whenever \(\nu_0 = \nu_1\). Indeed, \(W_{2,K}\) defines a metric since \(K\) is positive definite. Conversely, the contribution of the radial parts via \(W_{2,K}^2(\eta_{\nu_0},\eta_{\nu_1})\) is positive whenever the radial laws differ. In Proposition \ref{prop:spherical-reduction-intro} below, we show that \(W_{2,K}^2(\eta_{\nu_0},\eta_{\nu_1})\) coincides after spectral reduction with the optimal transport cost \(\cW_{2,\bm\kappa}^2(\beta_0,\beta_1)\). Consequently, we obtain
\begin{equation}\label{eq:Gamma-def}
 \Gamma(P_0,P_1)
 = W_{2,K}^2(\eta_{\nu_0},\eta_{\nu_1})
\end{equation}
for the generator effect defined in \eqref{eq:generator-effect-intro}.
\item  
   If $\widetilde R$ is a non-standardized radial
   variable with $0<\E[\widetilde R^2]<\infty$,
   set \(s=\sqrt{\E[\widetilde R^2]/d}\), \(R=\widetilde R/s,\) \(A_*=sA\).
   Then $m+\widetilde R AU=m+RA_*U$,
   with $\E[R^2]=d$ and covariance matrix
   $A_*A_*^{\trans}$.
\item
For fixed radial laws $\nu_0,\nu_1$, consider 
\begin{align}\label{def_GLi}
    G(K):=W_{2,K}^2(\eta_{\nu_0},\eta_{\nu_1}).
\end{align}
It is not difficult to see that the function $G$ is positively homogeneous, concave, and increasing with respect to the Loewner order. The latter means that, if \(\tilde{K}-K\) is positive semidefinite, then \(G(K)\leq G(\tilde K)\).
\end{enumerate}
\end{remark}

\subsection{The spectral reduction: eigenspaces and Dirichlet block radii}

The anisotropy in \eqref{eq:effect-equivalent} is completely encoded by the spectrum of $K$. To see this, let $\kappa_1,\ldots,\kappa_k$ be the distinct eigenvalues of $K$, let $E_j$ be the corresponding
mutually orthogonal eigenspaces, let $\pi_j$ denote the orthogonal projection onto $E_j$, and put
$n_j:=\dim E_j$. Then
\begin{equation}\label{eq:K-spectral-intro}
 K=\sum_{j=1}^k\kappa_j\pi_j,
 \qquad n_1+\cdots+n_k=d,
\end{equation}
and therefore the weighted spherical transport problem has cost 
\begin{align*}
 (z-w)^{\trans}K(z-w)
&=\sum_{j=1}^k \kappa_j
(z-w)^{\trans}\pi_j(z-w)  \\
&=\sum_{j=1}^k \kappa_j
(z-w)^{\trans}\pi_j^\trans\pi_j(z-w)
=\sum_{j=1}^k \kappa_j
\norm{\pi_jz-\pi_jw}^2,
\end{align*}
where we use for the second equality that \(\pi_j\) is an orthogonal projection.
The above identity explains why the eigenspaces are the correct objects to consider. On each eigenspace $E_j$, the matrix $K$ acts as multiplication by the corresponding eigenvalue $\kappa_j$, that is, $Kv=\kappa_j v$ for all $v\in E_j$.
Hence all directions within \(E_j\) carry the same quadratic weight, and simultaneous orthogonal transformations of \(\pi_jz\) and \(\pi_jw\) leave the corresponding cost unchanged. Directions belonging to different eigenspaces, on the other hand, must remain distinguishable because they are weighted differently.

Now, with each eigenspace \(E_j\), we consider the following block lengths, that we also denote as \emph{block radii}:
\[
r_j:=\|\pi_jz\|,
\qquad
s_j:=\|\pi_jw\|.
\]
Then, the Cauchy--Schwarz inequality gives
\begin{equation}\label{eq:block-reverse-triangle-intro}
\begin{aligned}
\|\pi_jz-\pi_jw\|^2
&=r_j^2+s_j^2-2\langle \pi_jz,\pi_jw\rangle
\geq r_j^2+s_j^2-2r_js_j
= (\|\pi_j z\| - \|\pi_j w\|)^2.
\end{aligned}
\end{equation}
Equality holds if and only if \(\pi_jz\) and \(\pi_jw\) are positively collinear, that is, whenever both are nonzero, they point in the same direction. 
Thus, for fixed block radii, the transport cost within each eigenspace is minimized by aligning the two vectors in the same direction. Consequently, the remaining optimization concerns only the joint coupling of the block radii.

Consequently, every coupling of \(Z\) and \(W\) induces a coupling of the \emph{block-radius vectors}
\begin{align}\label{def:p_K}
p_K(Z_0):=\bigl(\|\pi_1Z_0\|,\ldots,\|\pi_kZ_0\|\bigr),
\qquad
p_K(Z_1):=\bigl(\|\pi_1Z_1\|,\ldots,\|\pi_kZ_1\|\bigr),
\end{align}
and thus \eqref{eq:block-reverse-triangle-intro} yields the lower bound
\begin{align}\label{eq_lowerbproj}
    (Z_0-Z_1)^{\trans}K(Z_0-Z_1)
\geq
\sum_{j=1}^k \kappa_j
\bigl(\|\pi_jZ_0\|-\|\pi_jZ_1\|\bigr)^2.
\end{align}
Now, the crucial point is that the lower bound in \eqref{eq_lowerbproj} is attainable. 
The idea is to first choose for the two block-radius vectors \(p_K(Z_0)\) and \(p_K(Z_1)\) an optimal coupling, whose existence is guaranteed by the general optimal transport theory although its dependence structure is generally not known explicitly. 
We then choose independent uniform directions \(V_j\) in the eigenspaces \(E_j\), independently of the coupled block radii, using the same direction for both vectors within each eigenspace.

To make this construction precise, define
\begin{align}\label{def:beta_i}
    \beta_i
:=
\law\!\left(
R_i\sqrt{D_1},\ldots,R_i\sqrt{D_k}
\right), \qquad \text{for } (D_1,\ldots,D_k)
\sim
\mathrm{Dirichlet}\!\left(
\frac{n_1}{2},\ldots,\frac{n_k}{2}
\right)
\end{align}
independent of \(R_i\sim\nu_i\).
Let $(r,s)=((r_1,\ldots,r_k),(s_1,\ldots,s_k))\sim\gamma^*$ be an optimal coupling of
$\beta_0$ and $\beta_1$ for the weighted quadratic
cost $\sum_{j=1}^k\kappa_j(r_j-s_j)^2$.
Independently of $(r,s)$, choose independent random
vectors $V_1,\ldots,V_k$, where $V_j$ is uniformly
distributed on the unit sphere of $E_j$, and define
\begin{align}\label{def_Zast}
Z_0^*:=\sum_{j=1}^k r_jV_j ,
\qquad
Z_1^*:=\sum_{j=1}^k s_jV_j .
\end{align}
As a consequence of the spherical block decomposition (see Lemma \ref{lem:spherical-block} for details), we have
$Z_0^*\sim\eta_{\nu_0}$ and $Z_1^*\sim\eta_{\nu_1}$.
Moreover, since both vectors share the common direction $V_j$ within each eigenspace, we obtain
\begin{align*}
\|\pi_jZ_0^*-\pi_jZ_1^*\|^2
=\|(r_j-s_j)V_j\|^2
=(r_j-s_j)^2,
\end{align*}
where the last equality follows from \(\|V_j\| = 1\).
Consequently,
\begin{align*}
\mathbb E[(Z_0^*-Z^*_1)^\top K(Z_0^*-Z_1^*)]
=
\mathbb E_{\gamma^*}
\left[\sum_{j=1}^k\kappa_j(r_j-s_j)^2\right].
\end{align*}
Thus the lower bound in \eqref{eq_lowerbproj}
is attained, and the remaining optimization
concerns only the joint coupling of the
block-radius vectors.
The preceding arguments yield the following spectral reduction.

\begin{proposition}[Spectral reduction and weighted transport between spherical laws]
\label{prop:spherical-reduction-intro}
Let \(K\) be the Bures matrix in \eqref{eq:K-main-intro} with distinct eigenvalues \(\bm\kappa=(\kappa_1,\ldots,\kappa_k)\) and corresponding eigenspace dimensions \(n_1,\ldots,n_k\). For \(i\in \{0,1\}\), consider the distributions \(\beta_i\) in \eqref{def:beta_i}.
Then, we have
\begin{align}\label{eq:spherical-spectral-reduction-intro}
W_{2,K}^2(\eta_{\nu_0},\eta_{\nu_1})
=
\cW_{2,\bm\kappa}^2(\beta_0,\beta_1)
:=
\inf_{\gamma\in\Pi(\beta_0,\beta_1)}
\int_{\R_+^k\times\R_+^k}
\sum_{j=1}^k
\kappa_j(r_j-s_j)^2
\,\gamma(\dd r,\dd s).
\end{align}
\end{proposition}

\begin{remark}
\begin{enumerate}[label = (\alph*)]
    \item While \(\eta_{\nu_i}\) are spherical distributions on
\(\R^d\), the \(\beta_i\) in \eqref{def:beta_i} are distributions on the positive orthant \(\R_+^k\). They describe the joint distribution of the block radii of the spherical vectors with respect to the eigenspace decomposition of \(K\). In particular, the \(\beta_i\) are generally not spherical
distributions.
Hence, Proposition \ref{prop:spherical-reduction-intro} reduces a \(d\)-dimensional anisotropic quadratic transport problem between the two spherical laws to a \(k\)-dimensional transport problem on the positive orthant. Recall that \(k\leq d\) is the number of distinct eigenvalues of \(K\). The reduction removes all angular variables within the eigenspaces; what remains is an optimal transport problem between the corresponding block radii.
    \item The function \(G(K)\) in \eqref{def_GLi} is orthogonally invariant and concave, and hence Schur-concave in the eigenvalues of $K$;
see \cite[Proposition 3.C.2]{MarshallOlkinArnold2011}. 
The latter implies for positive definite matrices
$K$ and $\widetilde K$ with
$\tr K=\tr\widetilde K$, that
\begin{align}
    \sum_{j=1}^q\lambda_j^\downarrow(K)
\leq
\sum_{j=1}^q\lambda_j^\downarrow(\tilde{K}),
\quad q=1,\ldots,d-1,
\quad\Longrightarrow\quad
W_{2,K}^2(\eta_{\nu_0},\eta_{\nu_1})\geq W_{2,\tilde K}^2(\eta_{\nu_0},\eta_{\nu_1}),
\end{align}
where \(\lambda_j^\downarrow(K)\) are the eigenvalues of \(K\) in decreasing order, similarly for \(\tilde{K}\).
Consequently, for fixed $\tr K$, increasing spectral
dispersion cannot increase the generator effect.
In particular, the isotropic case
$K=(\tr K/d)I_d$ maximizes $G$ and attains
the Wong upper bound.
\end{enumerate}
\end{remark}

\section{Main results}\label{sec_main}

In this section, we present our main results.
First, we establish an exact decomposition of the
quadratic Wasserstein distance between elliptical
distributions in terms of the Bures--Wasserstein
distance and a lower-dimensional transport problem
between their block-radius distributions.
We also construct optimal couplings by lifting
optimal solutions of the reduced problem.

We then consider the special case where the
Bures matrix $K$ has only one eigenvalue, for which we obtain an explicit
transport cost and an optimal Monge map. Finally, we derive improved
lower and upper bounds for the general quadratic transport problem,
exploiting nested spectral projections and
conditionally comonotone couplings of the block radius vectors, respectively.

\subsection{Representation of $W_2^2(P_0,P_1)$ via spectral reduction}

We now combine the transport cost decomposition in Proposition \ref{lem:effect-intro} with
the spectral reduction established in
Proposition~\ref{prop:spherical-reduction-intro}
to obtain an exact representation of the quadratic
Wasserstein distance between elliptical distributions.
Moreover, we show that every optimal coupling of the
block-radius distributions can be lifted to an optimal
coupling of the original elliptical distributions
by using common uniform directions within each eigenspace.

\begin{theorem}[Bures decomposition and spectral reduction]\label{thm:main}
For \(i=0,1\), let $P_i=\cE_2(m_i,\Sigma_i,\nu_i)$,  where $\Sigma_i\in\Spd^d$ and $\nu_i\in \cR_d$. Let $K$ be the Bures matrix with eigenvalues \(\bm\kappa = (\kappa_1,\ldots,\kappa_k)\) and let \(\beta_i\) be given by \eqref{def:beta_i}. Then
\begin{equation}\label{eq:main-formula}
W_2^2(P_0,P_1)
=
\norm{m_0-m_1}^2
+d_{\BW}^2(\Sigma_0,\Sigma_1)
+\cW_{2,\bm\kappa}^2(\beta_0,\beta_1).
\end{equation}
Moreover, every optimal coupling $\gamma^*\in\Pi(\beta_0,\beta_1)$ can be lifted to an optimal coupling of $P_0$ and $P_1$ as follows. For \(Z_0^*\) and \(Z_1^*\) in \eqref{def_Zast}, define
\begin{equation}\label{eq:lift-ZW-intro}
X_0:=m_0+\Sigma_0^{1/2}Z_0^*,
\quad
X_1:=m_1+\Sigma_0^{-1/2}K Z_1^*.
\end{equation}
Then $X_0\sim P_0$, $X_1\sim P_1$, and $(X_0,X_1)$ is an optimal coupling in quadratic Wasserstein distance.
\end{theorem}

\begin{remark}
\label{rem:intro-generator-effect}
\begin{enumerate}[label = (\alph*)]
    \item The last term in \eqref{eq:main-formula} is the generator effect
\(\Gamma(P_0,P_1)\) defined in \eqref{eq:generator-effect-intro}. The block-radius distributions $\beta_i$ capture
both the global radial magnitude and its allocation
across the distinct eigenspaces of $K$.
Thus, the generator effect reflects the interaction
between the radial laws and the spectral geometry of $K$.
\item The number \(k\) of distinct eigenvalues of \(K\) determines the dimension of
the reduced transport problem. If \(K=\lambda I_d\), then \(k=1\) and only the
global radius remains. If \(K\) has simple spectrum, then \(k=d\) and the block
radii are the absolute coordinates in an eigenbasis; the signs can be matched
optimally, while the \(d\) absolute coordinates must still be coupled jointly.
\item The reduced problem \(\cW_{2,\bm\kappa}^2(\beta_0,\beta_1)\) does \emph{not} decompose into \(k\) separate one-dimensional transport problems. Indeed, the block radii satisfy \(\sum_{j=1}^k \bigl(R_i\sqrt{D_j}\bigr)^2=R_i^2\) and share the same Dirichlet vector \((D_1,\ldots,D_k)\). Hence their dependence
must be preserved, and coordinatewise optimal radial couplings need not combine
to a coupling with the prescribed marginal laws \(\beta_i\). In general, no
closed-form solution is available for this \(k\)-dimensional transport problem. Thus, Theorem~\ref{thm:main} provides an exact structural reduction rather
than a closed formula for the generator effect.
\end{enumerate}
\end{remark}

In the isotropic case, where the Bures matrix \(K\) is a multiple of the identity matrix \(I_d\), we can determine the optimal Monge map and its cost via Theorems \ref{lem:quadratic-optimality} and \ref{thm:main} as follows. We refer to Figure \ref{fig:optimal-monge-transport} for an illustration.

\begin{corollary}[Optimal Monge map]\label{cor:scalar-intro}
Assume that $K=\lambda I_d$ for some \(\lambda>0\).
Then
\begin{enumerate}[label = (\roman*)]
    \item \label{cor:scalar-intro1} \(W_2^2(P_0,P_1)=\norm{m_0-m_1}^2+d_{\BW}^2(\Sigma_0,\Sigma_1)
                       +\lambda W_2^2(\nu_0,\nu_1).\)
    \item \label{cor:scalar-intro2} If $\nu_0$ is non-atomic, an optimal Monge map \(T\colon \R^d\to \R^d\) is given by
\begin{equation}\label{eq:Monge-special-intro}
 T(x)=m_1+\lambda\frac{\tau(t(x))}{t(x)}\Sigma_0^{-1}(x-m_0),
 \qquad t(x)=\norm{\Sigma_0^{-1/2}(x-m_0)},
\end{equation}
\(P_0\)-almost surely,
where \(\tau=F_{\nu_1}^{-1}\circ F_{\nu_0}\).
\item \label{cor:scalar-intro3} 
\(T\) is the gradient of the proper convex potential
\begin{equation}\label{eq:potential-special-intro}
 \Phi(x)=m_1^{\trans}x+\lambda\int_0^{t(x)}\tau(s)\dd s,
\end{equation}
for \(P_0\)-almost all \(x\in \R^d\),
with the convention \(F_{\nu_1}^{-1}(0):=\inf\operatorname{supp}(\nu_1)\) and \(
F_{\nu_1}^{-1}(1):=\sup\operatorname{supp}(\nu_1)\).
\end{enumerate}
\end{corollary}

\begin{remark}
The explicit optimal Monge map in Corollary~\ref{cor:scalar-intro} requires
$K=\lambda I_d$ for some $\lambda>0$, or equivalently,
\[
\Sigma_1=\lambda^2\Sigma_0^{-1}.
\]
Thus, the covariance matrices need not be isotropic
or identical. Conversely, if $\Sigma_0=\Sigma_1=\Sigma$,
then $K=\Sigma$, which is isotropic only if
$\Sigma$ is a multiple of the identity.
Hence, the case of identical covariance matrices
is generally not covered by Corollary~\ref{cor:scalar-intro}.
Although an optimal Monge map exists by the convex gradient characterization in Theorem \ref{lem:quadratic-optimality}, no general explicit formula is available in this case.
\end{remark}

\subsection{Improved bounds for \(W_2(P_0,P_1)\)}

As discussed above, when the radial laws differ and $K$ has at least two distinct eigenvalues ($k\geq 2$), the generator effect
$\Gamma(P_0,P_1)=\cW_{2,\bm\kappa}^2(\beta_0,\beta_1)$
generally does not admit an explicit formula,
since the optimal coupling for the weighted
transport problem \eqref{eq:spherical-spectral-reduction-intro}
is not known in closed form.
This motivates the search for explicit and
computationally tractable lower and upper bounds
for the generator effect.

As a first consequence of Theorem~\ref{thm:main}, we obtain a simple radial-sensitive lower bound for the generator effect: It complements the elementary upper bound \(U_{\mathrm{Wong}}\) in \eqref{eq:bounded_Gamma} and improves the classical Gelbrich bound \(L_{\mathrm{Gelbrich}}= 0\) whenever \(\nu_0\ne \nu_1\) and \(k\geq 2\). Therefore, write
\begin{equation}\label{eq:delta-r-intro}
\delta_R^2
:=
W_2^2(\nu_0,\nu_1)
=
\int_0^1
\abs{F_{\nu_0}^{-1}(u)-F_{\nu_1}^{-1}(u)}^2\,\dd u
\end{equation}
for the quadratic Wasserstein distance between the radial laws. We denote by \(\lambda_{\min}(K)\) the smallest eigenvalue of \(K\) and by \(\tr K\) the trace of \(K\).

\begin{corollary}[Elementary bounds]\label{cor:intro-radial-bounds}
Under the assumptions of Theorem~\ref{thm:main}, the generator effect is bounded via 
\begin{align}\label{eq:bounded_Gamma}
    \lambda_{\min}(K)\,\delta_R^2
\leq
\Gamma(P_0,P_1)
\leq
\frac{\tr K}{d}\,\delta_R^2 = U_\mathrm{Wong}.
\end{align}
\end{corollary}

In the following, we substantially improve the elementary bounds
in \eqref{eq:bounded_Gamma}. First, we derive sharper lower bounds
by exploiting the spectral structure of the Bures  matrix \(K\). Then, we establish improved upper bounds by coupling the block radii conditionally comonotone.

Now, to determine improved lower bounds, let $q\in\{1,\ldots,d\}$.
Recall that the first component of a \(\mathrm{Dirichlet}(q/2,(d-q)/2)\)-distributed bivariate random vector is \(\mathrm{Beta}(q/2,(d-q)/2)\)-distributed; see Section \ref{sec:Dirichlet}. 
Hence, the Euclidean norm of the orthogonal projection
of a spherical random vector $RU$ onto any
$q$-dimensional subspace has the same distribution
as $R\sqrt{B_{q,d}}$, where
$B_{q,d}\sim\mathrm{Beta}(q/2,(d-q)/2)$ is
independent of $R$, with the convention $B_{d,d}=1$.
Now, denote the squared Wasserstein distance between two such projected spherical distributions by
\begin{equation}\label{eq:projection-deltas}
 \Delta_q:=W_2^2\bigl(\law(R_0\sqrt{B_{q,d}}),
                            \law(R_1\sqrt{B_{q,d}})\bigr),
 \qquad B_{q,d}\sim\mathrm{Beta}(q/2,(d-q)/2),
\end{equation}
where $B_{q,d}$ is independent of the corresponding radius \(R_i\). We set $\Delta_0=0$ and, for $q=d$, $B_{d,d}=1$. In particular, we have $\Delta_d=\delta_R^2$. 

For the improved lower bound, the key idea is to consider nested projections corresponding to increasing eigenvalues. Therefore, we order the distinct eigenvalues as
\begin{align}\label{eq:ord_kappa}
    0<\kappa_1<\cdots<\kappa_k
\end{align}
and consider the projections \(Q_j:= \pi_j+\cdots+\pi_k\) onto the nested subspaces \(E_j\oplus\cdots\oplus E_k\) with associated dimensions
\begin{align}\label{def:dim}
    d_j:=n_j+\cdots+n_k \qquad j=1,\ldots,k.
\end{align}
Using a telescoping sum argument, the matrix \(K\) can then be decomposed into
\begin{align}\label{eq:spectral-layer}
 K= \sum_{\ell=1}^k \kappa_\ell\pi_\ell = \kappa_1 \sum_{\ell=1}^k \pi_\ell + \sum_{j=2}^k (\kappa_j - \kappa_{j-1}) \sum_{\ell = j}^k \pi_\ell = \kappa_1 I_d+\sum_{j=2}^k(\kappa_j-\kappa_{j-1})Q_j.
\end{align}
The first term
\(\kappa_1 I_d\) acts on the whole space and yields the minimum eigenvalue bound in \eqref{eq:bounded_Gamma}, that is, \(\kappa_1 \Delta_d = \lambda_{\min}(K)\delta_R^2 \). The non-negative weights
\(\kappa_j-\kappa_{j-1}\) acting on the nested subspaces increase this lower bound as follows.

\begin{theorem}[Improved lower bound via nested projection]\label{thm:nested-bound}
Assume that the eigenvalues of \(K\) are ordered as in \eqref{eq:ord_kappa} and let \(d_j\) be given by \eqref{def:dim}, \(j\in \{1,\ldots,k\}\).
Then the generator effect \(\Gamma\) in \eqref{eq:Gamma-def} satisfies
\begin{equation}\label{eq:nested-bound}
  L_K:=\kappa_1\Delta_d+
       \sum_{j=2}^k(\kappa_j-\kappa_{j-1})\Delta_{d_j} \leq \Gamma(P_0,P_1).
\end{equation}
If $K$ has at least two distinct eigenvalues and $\nu_0\neq\nu_1$, the improvement over the lower bound
  \(\lambda_{\min}(K)\,\delta_R^2 = \kappa_1 \Delta_d\) is strict.
\end{theorem}

\begin{remark}[Comparison with separate-block bounds]
The nested-projection lower bound \(L_K\) dominates the
bound obtained by considering the eigenspaces
separately. Indeed, the additivity of the quadratic
cost across orthogonal subspaces implies
\begin{align*}
\Delta_{d_j}\geq
\sum_{\ell=j}^k\Delta_{n_\ell} \qquad \text{and thus} \qquad L_K\geq\sum_{j=1}^k\kappa_j\Delta_{n_j}.
\end{align*}
Consequently, nested projections yield a stronger bound
by accounting for the joint dependence between
block radii. The bound requires only $k$
one-dimensional Wasserstein distances and can therefore be evaluated efficiently. For $k=1$, we have \(L_K = \kappa_1 \Delta_d = \kappa_1 W_2^2(\nu_0,\nu_1) = \Gamma(P_0,P_1)\).
\end{remark}

To improve the elementary trace-based upper bound \(U_{\mathrm{Wong}}\) in \eqref{eq:bounded_Gamma}, we decompose \(\R^d\) into \(m \leq d\) pairwise orthogonal subspaces,
\begin{align}\label{dim_eigens}
    \R^d = E_1\oplus \cdots\oplus E_m.
\end{align}
Then we couple the corresponding block radii successively conditionally comonotone, following a Knothe--Rosenblatt-type construction. As our numerical study in Section \ref{sec:simulation} shows, this approach substantially improves
the upper bound \(U_{\mathrm{Wong}}\) even with a small number \(m\) of subspaces. 
In general, our improved upper bound \(U_K\) in Eq. \eqref{def_UkB} below can be
expressed as an $m$-dimensional integral, whose
numerical evaluation may become challenging for
large $m$. However, for certain radial families,
including multivariate Student-t distributions,
the bound reduces to a combination of several
one-dimensional integrals, regardless of the
ambient dimension $d$ or the number of groups $m$.

Now, to construct the improved upper bound, we denote the eigenvalues of \(K\), counted with multiplicity and sorted in decreasing order, by \(\lambda_1^\downarrow\geq \cdots \geq \lambda_d^\downarrow > 0\). We denote the sum of the \(q\) largest eigenvalues by
\begin{align}\label{sum_eigenvalue}
    \Lambda_q := \sum_{j=1}^q \lambda_j^\downarrow, \qquad q\in \{1,\ldots,d\}, \qquad  \Lambda_0 := 0.
\end{align}
Let $\mathbf b=(b_1,\ldots,b_m)$ be an ordered composition of $d$ which we interpret as the dimensions of the orthogonal subspaces in \eqref{dim_eigens}, i.e., $b_h = \dim(E_h)\in \{1,2,\ldots\}$ and $\sum_{h=1}^m b_h=d$. Similar to \eqref{def:beta_i}, consider an \(m\)-dimensional random vector
\[
(D_1,\ldots,D_m)
\sim\mathrm{Dirichlet}
\left(\frac{b_1}{2},\ldots,\frac{b_m}{2}\right),
\]
with $D_1=1$ if $m=1$.
For radial random variables $R_i\sim\nu_i$ independent of \((D_1,\ldots,D_m)\), consider the random vectors of block radii defined by
\begin{align}\label{eq_Zi_blockradii}
    Z_i=(Z_{i,1},\ldots,Z_{i,m}) := (R_i\sqrt{D_1},\ldots,R_i\sqrt{D_m}),
\qquad i=0,1.
\end{align}
Then we define the \emph{Knothe-Rosenblatt coupling} of \(Z^0\) and \(Z^1\) as the joint distribution of \(Z_0^{\mathrm{KR}}\) and \(Z_1^{\mathrm{KR}}\) recursively by
\begin{align}\label{constr_KR}
\begin{split}
    Z_{i,1}^{\mathrm{KR}} &:= F_{Z_{i,1}}^{-1}(U_1), \\
    Z_{i,h}^{\mathrm{KR}} &:= F_{Z_{i,h}|Z_{i,h-1}=Z_{i,h-1}^{\mathrm{KR}},\ldots,Z_{i,1}=Z_{i,1}^\mathrm{KR}}^{-1}(U_h), \quad h=2,\ldots,m,
\end{split}
\end{align}
for \(i\in \{0,1\}\), where \(U_1,\ldots,U_m\sim\mathcal{U}(0,1)\) are independent. Note that, by construction, \(Z_{0,h}^{\mathrm{KR}}\) and \(Z_{1,h}^{\rm{KR}}\) are \(m\)-dimensional random vectors that are conditionally comonotone given \(U_1,\ldots,U_{h-1}\).
We denote the quadratic cost between the \(h\)th components of the Knothe-Rosenblatt coupling by
\begin{align}\label{cost_C_hKR}
    C_h^{\mathbf b} := \E(Z_{0,h}^{\mathrm{KR}}-Z_{1,h}^{\mathrm{KR}})^2.
\end{align}
Weighting this cost with an average eigenvalue and
summing over all subspaces, we define the quantity
\begin{align}\label{def_UkB}
    U_K := \min_{m\in \{1,\ldots, d\} \atop {\mathbf{b} = (b_1,\ldots,b_m)\in \N^m \atop b_1+\cdots + b_m = d }}
    U_{K,\mathbf{b}}, \qquad U_{K,\mathbf{b}}:= \sum_{h=1}^m \frac{\Lambda_{a_h+b_h}-\Lambda_{a_h}}{b_h} \, C_h^{\mathbf{b}}, \qquad a_h:= \sum_{\ell=1}^{h-1} b_\ell.
\end{align}
The following theorem shows that \(U_K\) serves as an improved upper bound for the generator effect.

\begin{theorem}[Improved upper bound via Knothe-Rosenblatt coupling]
\label{thm:multi-conditional-upper}
Let $P_i=\cE_2(m_i,\Sigma_i,\nu_i)$, $i=0,1$,
with $\nu_0,\nu_1\in\cR_d$, and let $K$ be the
corresponding Bures matrix. With \(\Lambda_q\) in \eqref{sum_eigenvalue} denoting the sum of the \(q\) largest eigenvalues of \(K\), the generator effect satisfies
\begin{align}\label{thm:multi-conditional-upper1}
    \Gamma(P_0,P_1) \leq U_K
\leq
\frac{\tr K}{d}\delta_R^2.
\end{align}
\end{theorem}

\begin{remark}
In general, evaluating $U_{K,\mathbf b}$ for an
ordered composition $\mathbf b=(b_1,\ldots,b_m)$
requires multidimensional integration involving
up to $m$ variables. This is often considerably simpler than solving numerically
the $k$-dimensional optimal transport problem
in \eqref{eq:spherical-spectral-reduction-intro}. 
Further, computing $U_K$ generally involves considering
all $2^{d-1}$ ordered compositions of $d$.
Restricting the minimization to suitably
selected compositions can nevertheless yield
substantially improved upper bounds.
\end{remark}

\section{Fast computation of $U_K$ via dynamic programming}\label{sec:dyn}

In this section, we develop a computationally
efficient method for evaluating the conditionally
comonotone upper bound $U_K$ in
\eqref{def_UkB}.
Therefore, recall that its construction is based on the
Knothe--Rosenblatt coupling $(Z_0^{\mathrm{KR}},Z_1^{\mathrm{KR}})$
in \eqref{constr_KR} and the minimization of
$U_{K,\mathbf b}$ over all ordered compositions
$\mathbf b$ of $d$, corresponding to consecutive
blocks of the eigenvalues of $K$ arranged in
decreasing order.

For certain spherical distributions, the
Knothe--Rosenblatt coupling admits the product
representation
\begin{align}\label{eq_condKRa}
Z_{i,h}^{\mathrm{KR}}
=
\sqrt{c_i}\,g_{i;a_h,b_h}(U_h)
\prod_{\ell<h}
\sqrt{1+\varepsilon_i
g_{i;a_\ell,b_\ell}(U_\ell)^2},
\qquad i=0,1,
\end{align}
as established in Eq. \eqref{rep_produ_KR} in the proof of
Theorem~\ref{lem:KR-student-one-dimensional}
below. This representation allows us to reduce
the computation of $U_K$ to at most $d^2$
one-dimensional integrals. Moreover, dynamic
programming reduces the optimization over all
$2^{d-1}$ ordered compositions of $d$ to
$d(d+1)/2$ transitions, yielding an optimization algorithm with complexity $O(d^2)$.

We establish this representation for Gaussian and
Student-$t$ distributions, as well as for uniform
distributions on balls and spheres, including arbitrary
pairings of these four families. 
Recall the block-radius vectors \(Z_i=(Z_{i,1},\ldots,Z_{i,m})\) in
\eqref{eq_Zi_blockradii}.
The functions \(g_{i,a_h,b_h} = F_{Q_{i,a_h,b_h}}^{-1}\) in \eqref{eq_condKRa} are the standardized conditional quantile functions of the \(h\)th block radius \(Z_{i,h}\) given \(Z_{i,1},\ldots,Z_{i,h-1}\); see Lemma \ref{lem:KR-conditional-block-radii} for details. The distributions \(Q_{i;a_h,b_h}\) as well as the constants \(c_i\) and \(\varepsilon_i\) for the respective distributions are displayed in Table \ref{tab:KR-families}.

\begin{table}[t!]
\centering
\renewcommand{\arraystretch}{1.6}
\begin{tabular}{@{}l|ccc@{}}
\toprule
Spherical distribution \(\eta_{\nu}\) on \(\R^d\) (\(\nu\in \cR_d\)) & $c$ & $\varepsilon$
& $Q_{a,b}$\\
\midrule

Gaussian \(N(0,I_d)\)
& $1$ & $0$
& $\chi_b$ \\ \addlinespace[4pt]

\shortstack[l]{Student-$t_{\theta}(0,I_d),$ covariance-standardized,\\
$\theta>2$ degrees of freedom}

& $\theta-2$ & $1$
& $\displaystyle
\sqrt{\frac{\chi^2_b}{\chi^2_{\theta+a}}}$ \\ \addlinespace[4pt]

Uniform on a ball with radius \(\sqrt{d+2}\)
& $d+2$ & $-1$
& $\displaystyle
\sqrt{\mathrm{Beta}\left(
\frac b2,\frac{d-a-b}{2}+1
\right)}$ \\\addlinespace[4pt]

Uniform on a sphere with radius \(\sqrt{d}\)
& $d$ & $-1$
& $\displaystyle
\sqrt{\mathrm{Beta}\left(
\frac b2,\frac{d-a-b}{2}
\right)} $\\
\bottomrule
\end{tabular}
\caption{Constants and conditional distributions \(Q_{a_h,b_h} = \law(Z_h/\alpha_{h} \mid Z_{1:h-1}=r)\) for the block radius vector in \eqref{eq_Zi_blockradii} with normalization \(\alpha_{h}\) in \eqref{eq:KR-conditional-radius}. For the spherical distributions above, these conditional distributions are independent of \(r\).}
\label{tab:KR-families}
\end{table}

\begin{theorem}[Fast computation of $U_K$]
\label{lem:KR-student-one-dimensional}
For $i=0,1$, let
$P_i=\cE_2(0,\Sigma_i,\nu_i)$,
where $\Sigma_i$ is positive definite and \(\eta_{\nu_i}\) is a spherical distribution that belongs to one of the four families in Table \ref{tab:KR-families}.
For integers $a\geq0$ and $b\geq1$ with
$a+b\leq d$, define the one-dimensional integrals
\begin{align*}
I_{a,b}
&:=\int_0^1
g_{0;a,b}(u)g_{1;a,b}(u)\,\mathrm du, \qquad
J_{a,b}
:=\int_0^1
\sqrt{1+\varepsilon_0\,g_{0;a,b}(u)^2}\,
\sqrt{1+\varepsilon_1\,g_{1;a,b}(u)^2}
\,\mathrm du.
\end{align*}
where \(\varepsilon_i\) are constants and \((0,1) \ni u \mapsto g_{i;a,b}(u):=F_{Q_{i;a,b}}^{-1}(u)\) the generalized inverse distribution functions of the conditional block radial distributions \(Q_{i;a,b}\) specified in Table \ref{tab:KR-families}.
For $w_{a,b}:=(\Lambda_{a+b}-\Lambda_a)/b$,
define recursively
\begin{align}\label{eq:student-KR-DP}
B_d:=0 \qquad \text{and} \qquad
B_a:=
\max_{1\leq b\leq d-a}
\left\{
w_{a,b}I_{a,b}+J_{a,b}B_{a+b}
\right\},
\quad a=d-1,\ldots,0.
\end{align}
Then, we have
\begin{align*}
    U_K=2\tr K-2\sqrt{c_0c_1}\,B_0,
\end{align*}
where the respective constants \(c_i\) are given in Table \ref{tab:KR-families}.
\end{theorem}

\begin{remark}
\label{rem:student-KR-computation}
\begin{enumerate}[label=(\alph*)]
\item 
The maximization in \eqref{eq:student-KR-DP} determines
the optimal dimension of each successive block.
In particular, the optimal blocks need not be
one-dimensional and may have different dimensions.
For example, if $K=\lambda I_d$, the single-block partition
$\mathbf b=(d)$ is optimal, whereas anisotropic spectra
may favour finer partitions.
\item 

By Theorem~\ref{lem:KR-student-one-dimensional},
the computation of $U_K$ requires at most $d^2$
one-dimensional integrals, namely $d(d+1)/2$
integrals $I_{a,b}$ and $d(d-1)/2$ integrals
$J_{a,b}$, followed by $O(d^2)$ arithmetic
operations for the dynamic program
\eqref{eq:student-KR-DP}.

\item
For fixed radial laws and dimension $d$, the
one-dimensional integrals $I_{a,b}$ and $J_{a,b}$
are independent of the covariance matrices.
Hence, they can be precomputed and reused for
different covariance matrices, requiring only
an update of the spectral weights $w_{a,b}$
and the dynamic programming recursion.
\end{enumerate}
\end{remark}

\section{Numerical experiments}\label{sec:simulation}

In this section, we compare the minimum-eigenvalue lower bound $L_{\min}$,
the nested-projection lower bound $L_K$, and the conditionally
comonotone upper bound $U_K$ with the elementary Gelbrich--Wong
bounds. 
For selected settings, we provide numerical approximations
of the generator effect $\Gamma(P_0,P_1)$ by solving discretized versions of
the reduced optimal transport problem
\eqref{eq:spherical-spectral-reduction-intro},
provided that the results exhibit sufficient stability
under grid refinement.

\subsection{Settings}

All reported numbers concern the generator effect $\Gamma(P_0,P_1)$, which depends only on the spectrum of the Bures matrix \(K\) and the radial laws \(\nu_0\) and \(\nu_1\). Hence, we may choose w.l.o.g. $m_0=m_1=0$ and $\Sigma_0=\Sigma_1=K$.
The corresponding Bures matrix is then $K$ and
$W_2^2(P_0,P_1)=\Gamma(P_0,P_1)$.
The spectral models for \(K\) are specified in Table~\ref{tab:spectra-simple}, with the common normalization
\begin{equation}\label{eq:numerical-spectral-scaling}
\tr K=d^2.
\end{equation}
Scaling \(K\) transfers to the transport costs and bounds, but changes neither the optimal block composition nor the improvement percentage.
Under \eqref{eq:numerical-spectral-scaling}, the elementary bounds are
\[
L_{\mathrm{Gelbrich}}=0,
\qquad
L_{\min}=\lambda_{\min}(K)\delta_R^2,
\qquad
U_{\mathrm{Wong}}=d\,\delta_R^2,
\qquad
\delta_R^2=W_2^2(\nu_0,\nu_1).
\]

\begin{table}[t!]
\centering\small
\renewcommand{\arraystretch}{1.3}
\begin{tabular}{@{}lll@{}}
\toprule
Bures matrix & Eigenvalues $\lambda_j^\downarrow(K)$ & Distinct levels $k$\\
\midrule
$K_2$ & $3d$ ($d/5$ times), $d/2$ ($4d/5$ times) & $2$\\
$K_5$ & $2d,\ 3d/2,\ d,\ d/3,\ d/6$, each $d/5$ times & $5$\\
$K_{\mathrm{ar}}$ & $\displaystyle\frac{2d(d-j+1)}{d+1}$, $j=1,\ldots,d$ & $d$\\
$K_{\mathrm{pow}}$ & $\displaystyle\frac{d^2j^{-2}}{\sum_{\ell=1}^d\ell^{-2}}$, $j=1,\ldots,d$ & $d$\\
\bottomrule
\end{tabular}
\caption{Spectral models with $\tr K=d^2$ and for $d=2$, we set $K_2=\operatorname{diag}(3,1)$.
The eigenvalues are arranged in decreasing order and counted with multiplicity.}
\label{tab:spectra-simple}
\end{table}
The models \(K_2\) and \(K_5\) in Table \ref{tab:spectra-simple} have \(2\) and \(5\) distinct eigenvalues, respectively. 
The arithmetic spectrum \(K_{\mathrm{ar}}\) decreases linearly, whereas the
power spectrum \(K_{\mathrm{pow}}\) assigns much larger weights to its leading directions.
In particular, the ambient dimension $d$ and the dimension $k$ of the
reduced transport problem can be very different.

For \(P_i\), we consider Gaussian distributions (\(N\)), Student-\(t\) distributions with three and five degrees of freedom (\(t_3\) and \(t_5\)), the uniform distribution on the ball of radius
$\sqrt{d+2}$ (Ball), and the uniform distribution on the sphere of
radius $\sqrt d$ (Sphere). Their radial variables can be represented as
\begin{align*}
R_N&=\sqrt{H_d},
&R_{t_\theta}&=\sqrt{\frac{(\theta-2)H_d}{V_\theta}},
\qquad \theta\in\{3,5\},\\
R_{\mathrm{Ball}}&=\sqrt{d+2}\,W^{1/d},
&R_{\mathrm{Sphere}}&=\sqrt d,
\end{align*}
where $H_d\sim\chi_d^2$, $V_\theta\sim\chi_\theta^2$, and
$W\sim\mathcal U(0,1)$ are independent.
Each radial law is in \(\cR_d\) and thus the associated radial variable \(R\) satisfies $\mathbb E[R^2]=d$.

Table~\ref{tab:intro-highlights} in Section \ref{sec:intro} presents the bounds for the generator effect between the normal and Student-$t_3$ distributions.
Table~\ref{tab:numerical-results} complements the introduction
with the pairs $N$--Sphere, Sphere--Ball, and $t_5$--$t_3$.
We consider small-, moderate-, and high-dimensional settings with dimensions \(d\in \{2,10,100,10^6\}\).

In both tables, ``Improvement'' quantifies the reduction
in the classical Gelbrich--Wong interval achieved by our
new bounds:
\begin{equation}\label{eq:simulation-GW-improvement-new}
\mathrm{Improvement}
:=
1-\frac{U_K-L_K}{U_{\mathrm{Wong}}-L_{\mathrm{Gelbrich}}}
=1-\frac{U_K-L_K}{d\,\delta_R^2}.
\end{equation}

\paragraph{Evaluation of the bounds.}

We compute the lower bound $L_K$ in \eqref{eq:nested-bound}
directly using one-dimensional quadrature.
For the upper bound $U_K$, we apply
Theorem~\ref{lem:KR-student-one-dimensional},
which expresses the required KR costs in terms of
one-dimensional integrals and allows optimization
over ordered block compositions by dynamic programming.

In high dimensions, we evaluate the projection costs
only for selected dimensions and restrict the optimization
of the upper bound to a subset of ordered block decompositions.
The following lemma allows us to control the resulting
approximation of the lower bound.

\begin{lemma}[Monotonicity of normalized projection costs]
\label{lem:normalized-projection-costs}
Let $\nu_0,\nu_1\in\cR_d$ and let $\Delta_q$,
$q=1,\ldots,d$, be the squared transport cost defined in
\eqref{eq:projection-deltas}. Then we have
\begin{equation}\label{eq:normalized-projection-costs}
\frac{\Delta_q}{q}
\leq
\frac{\Delta_r}{r}\qquad \text{for }
 1\leq q\leq r\leq d.
\end{equation}
\end{lemma}

For large $d$, we evaluate the projection costs
$\Delta_q$ only at selected dimensions
$1=q_1<\cdots<q_s=d$.
We then define
\begin{equation}\label{eq:simulation-bracket}
L^-:=
\omega_d\Delta_d+
\sum_{j=1}^{s-1}
\frac{\Delta_{q_j}}{q_j}
\sum_{q=q_j}^{q_{j+1}-1}q\omega_q,
\end{equation}
where
$\omega_q=\lambda_q^\downarrow(K)
-\lambda_{q+1}^\downarrow(K)$
and $\lambda_{d+1}^\downarrow(K)=0$.
Then, Lemma~\ref{lem:normalized-projection-costs} yields $L^-\le L_K$.

For a conservative approximation of the upper bound in high dimensions, we consider a subset 
$\mathcal A\subseteq\{0,\ldots,d\}$ that contains $0$ and $d$. Then, restricting the dynamic program in \eqref{eq:student-KR-DP} to transitions between indices in
$\mathcal A$ gives $U^+$.
Consequently,
\begin{equation}\label{eq:simulation-conservative-bounds}
L^-\le L_K\le\Gamma(P_0,P_1)\le U_K\le U^+
\le U_{\mathrm{Wong}}.
\end{equation}

For $d=10^6$, we approximate the lower bound using
$3663$ selected projection dimensions and restrict the
upper-bound optimization to $169$ block boundaries,
resulting in $14\,196$ possible transitions.
A dagger in the tables indicates $L^-$ or $U^+$.

\begin{table}[t]
\centering
\scriptsize
\setlength{\tabcolsep}{1.7pt}
\renewcommand{\arraystretch}{1.08}
\begin{tabular}{@{}rllrrrrrrr@{}}
\toprule
$d$ & Radial pair & \shortstack{Bures\\matrix} & $k$
& $L_{\min}$ & $L_K$ & $\Gamma(P_0,P_1)$
& $U_K$ & $U_{\mathrm{Wong}}$ & Improvement \\
\midrule
$2$ & $N$--Sphere & $K_2$ & $2$ & $0.4551$ & $0.6614$ & $\approx 0.77$ & $0.8727$ & $0.9102$ & $76.78\%$ \\
$2$ & $N$--Sphere & $K_{\mathrm{ar}}$ & $2$ & $0.6068$ & $0.7443$ & $\approx 0.85$ & $0.9102$ & $0.9102$ & $81.78\%$ \\
$2$ & $N$--Sphere & $K_{\mathrm{pow}}$ & $2$ & $0.3641$ & $0.6116$ & $\approx 0.72$ & $0.7807$ & $0.9102$ & $81.42\%$ \\
$2$ & Sphere--Ball & $K_2$ & $2$ & $0.2288$ & $0.2905$ & $\approx 0.36$ & $0.4227$ & $0.4575$ & $71.09\%$ \\
$2$ & Sphere--Ball & $K_{\mathrm{ar}}$ & $2$ & $0.305$ & $0.3461$ & $\approx 0.41$ & $0.4575$ & $0.4575$ & $75.66\%$ \\
$2$ & Sphere--Ball & $K_{\mathrm{pow}}$ & $2$ & $0.183$ & $0.257$ & $\approx 0.32$ & $0.3629$ & $0.4575$ & $76.87\%$ \\
$2$ & $t_{5}$--$t_{3}$ & $K_2$ & $2$ & $0.1398$ & $0.2716$ & -- & $0.2795$ & $0.2795$ & $97.16\%$ \\
$2$ & $t_{5}$--$t_{3}$ & $K_{\mathrm{ar}}$ & $2$ & $0.1863$ & $0.2742$ & -- & $0.2795$ & $0.2795$ & $98.11\%$ \\
$2$ & $t_{5}$--$t_{3}$ & $K_{\mathrm{pow}}$ & $2$ & $0.1118$ & $0.27$ & -- & $0.2772$ & $0.2795$ & $97.43\%$ \\
\addlinespace[7pt]
$10$ & $N$--Sphere & $K_2$ & $2$ & $2.465$ & $2.7411$ & $\approx 3.17$ & $3.2666$ & $4.93$ & $89.34\%$ \\
$10$ & $N$--Sphere & $K_5$ & $5$ & $0.8217$ & $1.9781$ & -- & $3.2228$ & $4.93$ & $74.75\%$ \\
$10$ & $N$--Sphere & $K_{\mathrm{ar}}$ & $10$ & $0.8964$ & $2.4005$ & -- & $4.3512$ & $4.93$ & $60.43\%$ \\
$10$ & Sphere--Ball & $K_2$ & $2$ & $0.4141$ & $0.4225$ & $\approx 0.50$ & $0.522$ & $0.8282$ & $87.98\%$ \\
$10$ & Sphere--Ball & $K_5$ & $5$ & $0.138$ & $0.1897$ & -- & $0.5064$ & $0.8282$ & $61.76\%$ \\
$10$ & Sphere--Ball & $K_{\mathrm{ar}}$ & $10$ & $0.1506$ & $0.2376$ & -- & $0.8282$ & $0.8282$ & $28.69\%$ \\
$10$ & $t_{5}$--$t_{3}$ & $K_2$ & $2$ & $3.9601$ & $7.454$ & -- & $7.7978$ & $7.9201$ & $95.66\%$ \\
$10$ & $t_{5}$--$t_{3}$ & $K_5$ & $5$ & $1.32$ & $7.6142$ & -- & $7.8022$ & $7.9201$ & $97.63\%$ \\
$10$ & $t_{5}$--$t_{3}$ & $K_{\mathrm{ar}}$ & $10$ & $1.44$ & $7.6856$ & -- & $7.8794$ & $7.9201$ & $97.55\%$ \\
\addlinespace[7pt]
$100$ & $N$--Sphere & $K_2$ & $2$ & $24.9684$ & $26.501$ & $\approx 30.9$ & $31.782$ & $49.9367$ & $89.42\%$ \\
$100$ & $N$--Sphere & $K_5$ & $5$ & $8.3228$ & $17.2394$ & -- & $30.3787$ & $49.9367$ & $73.69\%$ \\
$100$ & $N$--Sphere & $K_{\mathrm{ar}}$ & $100$ & $0.9888$ & $17.5728$ & -- & $35.285$ & $49.9367$ & $64.53\%$ \\
$100$ & Sphere--Ball & $K_2$ & $2$ & $0.4902$ & $0.4908$ & $\approx 0.59$ & $0.6129$ & $0.9803$ & $87.55\%$ \\
$100$ & Sphere--Ball & $K_5$ & $5$ & $0.1634$ & $0.1697$ & -- & $0.6101$ & $0.9803$ & $55.07\%$ \\
$100$ & Sphere--Ball & $K_{\mathrm{ar}}$ & $100$ & $0.0194$ & $0.0492$ & -- & $0.9515$ & $0.9803$ & $7.97\%$ \\
$100$ & $t_{5}$--$t_{3}$ & $K_2$ & $2$ & $419.9412$ & $827.9996$ & -- & $839.3492$ & $839.8823$ & $98.65\%$ \\
$100$ & $t_{5}$--$t_{3}$ & $K_5$ & $5$ & $139.9804$ & $833.67833$ & -- & $837.75444$ & $839.8823$ & $99.51\%$ \\
$100$ & $t_{5}$--$t_{3}$ & $K_{\mathrm{ar}}$ & $100$ & $16.6313$ & $834.19902$ & -- & $838.82535$ & $839.8823$ & $99.45\%$ \\
\addlinespace[7pt]
$10^6$ & $N$--Sphere & $K_5$ & $5$ & $83\,333.3229$ & $169\,315.7591$ & -- & ${301\,316.9529}^{\dagger}$ & $499\,999.9375$ & $73.60\%$ \\
$10^6$ & $N$--Sphere & $K_{\mathrm{ar}}$ & $10^6$ & $1$ & ${166\,397.7161}^{\dagger}$ & -- & ${328\,483.2623}^{\dagger}$ & $499\,999.9375$ & $67.58\%$ \\
$10^6$ & $N$--Sphere & $K_{\mathrm{pow}}$ & $10^6$ & $0.304$ & ${2.6728}^{\dagger}$ & -- & ${4.761}^{\dagger}$ & $499\,999.9375$ & $99.999582\%$ \\
$10^6$ & Sphere--Ball & $K_5$ & $5$ & $0.1667$ & $0.1667$ & -- & ${0.625}^{\dagger}$ & $1$ & $54.17\%$ \\
$10^6$ & Sphere--Ball & $K_{\mathrm{ar}}$ & $10^6$ & $2.0000\times10^{-6}$ & ${1.4168\times10^{-5}}^{\dagger}$ & -- & ${0.9635}^{\dagger}$ & $1$ & $3.65\%$ \\
$10^6$ & Sphere--Ball & $K_{\mathrm{pow}}$ & $10^6$ & $6.0793\times10^{-7}$ & ${6.0794\times10^{-7}}^{\dagger}$ & -- & ${2.4317\times10^{-6}}^{\dagger}$ & $1$ & $99.999818\%$ \\
$10^6$ & $t_{5}$--$t_{3}$ & $K_5$ & $5$ & $1.4110\times10^{10}$ & $8.46604524\times10^{10}$ & -- & ${8.46604984\times10^{10}}^{\dagger}$ & $8.4661\times10^{10}$ & $99.999946\%$ \\
$10^6$ & $t_{5}$--$t_{3}$ & $K_{\mathrm{ar}}$ & $10^6$ & $169\,320.8734$ & ${8.46604522\times10^{10}}^{\dagger}$ & -- & ${8.46605095\times10^{10}}^{\dagger}$ & $8.4661\times10^{10}$ & $99.999932\%$ \\
$10^6$ & $t_{5}$--$t_{3}$ & $K_{\mathrm{pow}}$ & $10^6$ & $51\,467.4567$ & ${7.0736\times10^{10}}^{\dagger}$ & -- & ${7.4540\times10^{10}}^{\dagger}$ & $8.4661\times10^{10}$ & $95.51\%$ \\
\bottomrule
\end{tabular}
\caption{Bounds and numerical OT approximations
for the selected radial pairs and spectra in
Table~\ref{tab:spectra-simple}.
Values in the $\Gamma(P_0,P_1)$ column are numerical estimates for the generator effect.
All matrices have trace $d^2$.
A dagger (\(\dagger\)) denotes $L^-$ or $U^+$ for $d=10^6$;
a dash (--) indicates that no sufficiently stable OT
approximation is reported.
Improvement is the reduction in the
Gelbrich--Wong interval width; see \eqref{eq:simulation-GW-improvement-new}.}
\label{tab:numerical-results}
\end{table}

\subsection{Discussion of the numerical results}

Table~\ref{tab:intro-highlights} illustrates the effectiveness
of our bounds for normal and Student-$t_3$ distributions.
The improvements are already substantial in dimension $2$
and exceed $99.95\%$ in dimension $10^6$ for $K_2$
and $K_{\mathrm{ar}}$.
In these examples, the nested-projection lower bound \(L_K\)
captures most of the radial discrepancy, while the
KR upper bound \(U_K\) provides further improvements,
particularly for anisotropic spectra.

Table~\ref{tab:numerical-results} extends the comparison
to other radial laws. Numerical approximations of
$\Gamma(P_0,P_1)$ are reported when they are 
 sufficiently stable under grid refinement. A dash indicates that no reliable
approximation was obtained; some numerical transport
couplings even yielded costs exceeding $U_K$.

\subsubsection*{Small dimensions}

For $d=2$, the bounds are particularly tight for
$t_5$--$t_3$, with improvements exceeding $97\%$.
For the other radial pairs, the reductions are less
pronounced. In several cases, $U_K=U_{\mathrm{Wong}}$,
so the improvement is entirely due to the lower bound.

At $d=10$, the spectral structure becomes more influential.
For Sphere--Ball, the improvement ranges from nearly
$88\%$ under $K_2$ to less than $29\%$ under
$K_{\mathrm{ar}}$. Thus, the effectiveness of the bounds
depends not only on the radial laws but also on
how their discrepancy interacts with the spectrum of $K$.

\subsubsection*{Moderate dimension}

At $d=100$, the bounds for $t_5$--$t_3$ become
even tighter, with improvements exceeding $98\%$
for all three spectra. By contrast, for Sphere--Ball
under $K_{\mathrm{ar}}$, the improvement drops below
$8\%$. Here, the projection-based lower bound remains
relatively weak, while the KR upper bound is close
to Wong's bound. Hence, increasing dimension
does not necessarily improve the tightness of our bounds.

\subsubsection*{High dimension}

At $d=10^6$, the contrasting trends become even more
pronounced. For $t_5$--$t_3$, the conservative improvements
exceed $99.9999\%$ under $K_5$ and $K_{\mathrm{ar}}$,
whereas Sphere--Ball under $K_{\mathrm{ar}}$ achieves
only $3.65\%$.
The strongly concentrated power spectrum, by contrast,
yields improvements above $99.9995\%$ for both
$N$--Sphere and Sphere--Ball.
This illustrates the strong influence of spectral
anisotropy on the quality of the bounds.
However, a nearly complete reduction of the
Gelbrich--Wong interval does not necessarily imply
a comparably small relative error in the generator effect.

\paragraph{Asymptotic scales.}
The different magnitudes of the transport costs can be
explained by the asymptotic behavior of the radial laws.
As $d\to\infty$, we have
\[
\delta_R^2 =
\begin{cases}
\Theta(d), & N\text{--}t_3,\quad t_5\text{--}t_3,\\
\Theta(1), & N\text{--Sphere},\quad N\text{--Ball},\\
\Theta(d^{-1}), & \text{Sphere--Ball}.
\end{cases}
\]
Since $\tr K=d^2$, the Wong bound satisfies
$U_{\mathrm{Wong}}=d\,\delta_R^2$ and is therefore
of order $d^2$, $d$, or $1$, respectively.
For $K_2$ and $K_5$, the smallest eigenvalue is
of order $d$, so the minimum-eigenvalue lower bound
shows that $\Gamma(P_0,P_1)$ has the same asymptotic
order. The same orders hold for the Student pairs
and $N$--Sphere under $K_{\mathrm{ar}}$.
For the strongly anisotropic power spectrum, however,
the transport costs can be substantially smaller,
as illustrated by the $N$--Sphere and Sphere--Ball
comparisons. Thus, the magnitude of the generator
effect depends on both the radial discrepancy and
the spectral structure of $K$.

\paragraph{Computation time.}

Including the quadrature checks, the computation took approximately $0.35$ seconds for $d=2$, $0.86$ seconds for $d=10$, $7.35$ seconds for $d=100$,
and $113$ seconds for the restricted calculation at $d=10^6$.
The R implementation used one worker for $d=2,10$ and eight
workers for the larger dimensions. 
The reported computation times include all five
   radial pairs, including $N$--Ball, which is omitted
   from the tables because of its similarity to
   $N$--Sphere.

\appendix

\section{Auxiliary results}\label{app_a}

\subsection{Elliptical distributions and their radial laws}

For general background on elliptical distributions, we refer to \cite{CambanisHuangSimons1981,FangKotzNg1990}. In this paper, we consider elliptical distributions with finite second moment. The normalization $\E[R^2]=d$ of the radial variable due to the class \(\cR_d\) in \eqref{eq:Rclass-intro} allows us to identify the matrix \(\Sigma\) as the covariance matrix of \(X = \mu+RAU\), \(AA^\top = \Sigma\), as follows.

\begin{lemma}[Normalization and uniqueness of the standardized radius]\label{lem:ell-normalization}
Let $X\sim\cE_2(m,\Sigma,\nu)$ with $\Sigma\in\Spd^d$ and $\nu\in\cR_d$.
Then
\begin{equation}\label{eq:mean-cov}
\E[X]=m,
\qquad
\Cov(X)=\Sigma, \qquad \text{and} \qquad \norm{\Sigma^{-1/2}(X-m)}\eqd R, \quad \text{for } R\sim\nu.
\end{equation}
The standardized radial law $\nu$ is uniquely determined by the distribution of $X$ together with $(m,\Sigma)$. Further, the law in \eqref{eq:elliptical-def-intro} is independent of the chosen factor $A$ satisfying $AA^{\trans}=\Sigma$.
\end{lemma}

\begin{proof}
Since $U$ is centered and \(\E[UU^{\trans}]=I_d/d\),
independence and $\E[R^2]=d$ give
\[
\E[ARU]=0,
\qquad
\Cov(ARU)
=A\E[R^2UU^{\trans}]A^{\trans}=AA^{\trans}=\Sigma.
\]
Let $A$ be any nonsingular factor with $AA^{\trans}=\Sigma$. Then \(
Q:=\Sigma^{-1/2}A\)
satisfies $QQ^{\trans}=I_d$ and is therefore orthogonal. Hence
\[
\norm{\Sigma^{-1/2}(X-m)}
=
\norm{QRU}
=R
\quad\text{a.s.},
\]
which proves \eqref{eq:mean-cov}. If $A_1A_1^{\trans}=A_2A_2^{\trans}=\Sigma$, then $A_2=A_1Q$ for an orthogonal matrix $Q$. Since $QU\eqd U$, the two factors generate the same law.
\end{proof}

\subsection{Dirichlet distribution}\label{sec:Dirichlet}

For parameters $\alpha_1,\ldots,\alpha_k>0$, a \(k\)-dimensional random vector
$D=(D_1,\ldots,D_k)$ has the \emph{Dirichlet distribution}
$\mathrm{Dirichlet}(\alpha_1,\ldots,\alpha_k)$ if it is supported on the simplex
\[
\Delta_{k-1}:=\left\{d\in[0,1]^k:\sum_{j=1}^k d_j=1\right\}
\]
and has density
\[
f_D(d_1,\ldots,d_k)
=
\frac{\Gamma(\alpha_1+\cdots+\alpha_k)}
     {\prod_{j=1}^k\Gamma(\alpha_j)}
\prod_{j=1}^k d_j^{\alpha_j-1}
\]
with respect to Lebesgue measure in the first
$k-1$ coordinates, where
$d_k=1-\sum_{j=1}^{k-1}d_j$.
Equivalently, if \(\{Y_j\}_{j=1,\ldots,k}\) are independent chi-squared random variables with
\(2\alpha_j\) degrees of freedom, that is,
\(
Y_j\sim\chi^2_{2\alpha_j}\), \(j=1,\ldots,k\)
then
\begin{align}\label{rep_Dirichlet}
    D_j=\frac{Y_j}{Y_1+\cdots+Y_k},
\qquad j=1,\ldots,k,
\end{align}
defines a random vector \(D=(D_1,\ldots,D_k)
\sim
\mathrm{Dirichlet}(\alpha_1,\ldots,\alpha_k)\); see e.g. \cite[Chapter~49]{KotzBalakrishnanJohnson2000}.
In particular, the marginal laws are given by
\[
D_j\sim
\mathrm{Beta}\!\left(\alpha_j,\sum_{\ell\neq j}\alpha_\ell\right).
\]
For $k=1$ we use the convention $D_1=1$ almost surely. Further, \(\mathrm{Beta}\!\left(\alpha,0\right)\) is the Dirac measure in \(1\). The representation in \eqref{rep_Dirichlet} will be
applied in Lemma~\ref{lem:spherical-block} with $\alpha_j=n_j/2$.
We also use that
\begin{align}\label{eq_betachi}
    \frac{X}{X+Y}\sim \mathrm{Beta}\left(\frac m2, \frac n 2\right) \qquad \text{and} \qquad \frac{X}{X+Y}\perp (X+Y)
\end{align}
for independent random variables \(X\sim \chi_m^2\) and \(Y\sim \chi_n^2\).

\subsection{Bures alignment and a Procrustes factor}

For two matrices $C,D\in\R^{d\times d}$, the orthogonal \emph{Procrustes problem} asks 
for an orthogonal matrix $Q\in O(d)$ that minimizes
\begin{align}\label{eq:procrustes-general}
    \norm{C-DQ}_{\mathrm F}^2.
\end{align}
Here, \(\norm{\cdot}_{\mathrm F}\) denotes the Frobenius norm, \(
\norm{A}_{\mathrm F}
:=
\left(\sum_{i,j=1}^d A_{ij}^2\right)^{1/2}
=
\sqrt{\tr(A^{\trans}A)}.
\)
This problem is particularly relevant for the optimal transport between the covariances of elliptical distributions because the matrix factors of the covariances are non-unique.
If \(DD^{\trans}=\Sigma\), then,
for every \(Q\in O(d)\),
\[
(DQ)(DQ)^{\trans}
=
DQQ^{\trans}D^{\trans}
=
DD^{\trans}
=
\Sigma.
\]
Thus right multiplication by an orthogonal matrix changes the factor
\(D\), but not the positive definite matrix generated by it. The
Procrustes problem selects, among all factors \(DQ\) of the same matrix,
the one that is closest to \(C\) in Frobenius norm.
Moreover,
\begin{align*}
\norm{C-DQ}_{\mathrm F}^2
&=
\tr\!\left((C-DQ)^{\trans}(C-DQ)\right)
=
\tr(CC^{\trans})
+\tr(DD^{\trans})
-2\tr(C^{\trans}DQ).
\end{align*}
Hence minimizing \eqref{eq:procrustes-general} is equivalent to maximizing
\(\tr(C^{\trans}DQ)\) over \(Q\in O(d)\).
For the covariance factors \(C=\Sigma_0^{1/2}\) and \(D=\Sigma_1^{1/2}\), 
the classical Procrustes identity gives
\begin{align*}
\max_{Q\in O(d)}
\tr\!\left(\Sigma_0^{1/2}\Sigma_1^{1/2}Q\right)
&=
\tr\!\left[
\left(
\Sigma_0^{1/2}\Sigma_1\Sigma_0^{1/2}
\right)^{1/2}
\right];
\end{align*}
see e.g. \cite[Theorem 1]{BhatiaJainLim2019}.
Consequently,
\begin{align}\label{eq:procrustes-trace2}
    \min_{Q\in O(d)}
\norm{\Sigma_0^{1/2}-\Sigma_1^{1/2}Q}_{\mathrm F}^2
=
d_{\BW}^2(\Sigma_0,\Sigma_1).
\end{align}
Thus the squared Bures--Wasserstein distance is the minimal squared
Frobenius distance between covariance factors after optimal orthogonal
alignment.

For \(\Sigma_0,\Sigma_1\in\Spd^d\), recall the Bures  matrix in \eqref{eq:K-main-intro} given by
\[
K=
\left(
\Sigma_0^{1/2}\Sigma_1\Sigma_0^{1/2}
\right)^{1/2}
\]

The following lemma provides optimally
aligned covariance factors, which play a key role
in the decomposition of the transport cost.

\begin{lemma}[Bures-aligned covariance factors]
\label{lem:bures-factors}
Let \(A:=\Sigma_0^{1/2}\) and \(B:=\Sigma_0^{-1/2}K\).
Then
\begin{equation}\label{eq:AB-properties}
AA^{\trans}=\Sigma_0,
\qquad
BB^{\trans}=\Sigma_1,
\qquad
A^{\trans}B=K.
\end{equation}
Moreover, the matrix \(Q^*
:=
\Sigma_1^{-1/2}\Sigma_0^{-1/2}K\)
is orthogonal and is an optimizer of the orthogonal Procrustes problem. In
particular,
\begin{align}\label{eq:B-procrustes}
B=\Sigma_1^{1/2}Q^* \qquad \text{and} \qquad \norm{A-B}_{\mathrm F}^2
=
d_{\BW}^2(\Sigma_0,\Sigma_1)
=
\min_{Q\in O(d)}
\norm{\Sigma_0^{1/2}-\Sigma_1^{1/2}Q}_{\mathrm F}^2.
\end{align}
\end{lemma}

\begin{proof}
The identities in \eqref{eq:AB-properties} are immediate. Moreover \(Q^*\) is orthogonal because \[
Q^*(Q^*)^{\trans}
=
\Sigma_1^{-1/2}\Sigma_0^{-1/2}
K^2
\Sigma_0^{-1/2}\Sigma_1^{-1/2}
=
I_d.
\]
Further, \(\Sigma_1^{1/2}Q^*
=
\Sigma_0^{-1/2}K
=
B,
\)
which proves the first part in \eqref{eq:B-procrustes}.
Finally,
\begin{align*}
\norm{A-B}_{\mathrm F}^2
&=
\tr(AA^{\trans})
+\tr(BB^{\trans})
-2\tr(A^{\trans}B)
=
\tr(\Sigma_0)
+\tr(\Sigma_1)
-2\tr(K)
=
d_{\BW}^2(\Sigma_0,\Sigma_1),
\end{align*}
where the last equality is by the definition of the Bures-Wasserstein distance in \eqref{eq:bures-def-intro} and \(K\).
The last equality in \eqref{eq:B-procrustes} follows from the Procrustes identity in \eqref{eq:procrustes-trace2}.
\end{proof}

\section{Proofs}

\subsection{Proofs of Section~\ref{sec_22}}\label{sec:proofs-intro}


\begin{proof}[Proof of Proposition~\ref{lem:effect-intro}]
By Lemma~\ref{lem:bures-factors}, we have \(AA^{\trans}=\Sigma_0,\) \(BB^{\trans}=\Sigma_1,\) and \(A^{\trans}B=K.\)
Since \(A\) and \(B\) are invertible, couplings of \(P_0\) and \(P_1\)
are in one-to-one correspondence with couplings of
\(\eta_{\nu_0}\) and \(\eta_{\nu_1}\) through
\[
X_0=m_0+AZ_0,
\qquad
X_1=m_1+BZ_1.
\]
For every such coupling, \eqref{eq:mean-separation-intro} and
\eqref{eq:couplingwise-preview-intro} give
\begin{align*}
\E\norm{X_0-X_1}^2
&=
\norm{m_0-m_1}^2
+d_{\BW}^2(\Sigma_0,\Sigma_1)
+\E\bigl[(Z_0-Z_1)^{\trans}K(Z_0-Z_1)\bigr].
\end{align*}
Taking the infimum over all couplings of \(X_0\) and \(X_1\) therefore yields
\[
W_2^2(P_0,P_1)
=
\norm{m_0-m_1}^2
+d_{\BW}^2(\Sigma_0,\Sigma_1)
+W_{2,K}^2(\eta_{\nu_0},\eta_{\nu_1}),
\]
which proves \eqref{eq:effect-formula}.

Finally, by the linear change of variables \(z\mapsto K^{1/2}z\),
\[
W_{2,K}^2(\eta_{\nu_0},\eta_{\nu_1})
=
W_2^2\bigl(
(K^{1/2})\push\eta_{\nu_0},
(K^{1/2})\push\eta_{\nu_1}
\bigr) = W_2^2(\cE_2(0,K,\nu_0),\cE_2(0,K,\nu_1)),
\]
where we use for the second equality \(K^{1/2}(K^{1/2})^{\trans}=K\) the definition of elliptical distributions in \eqref{eq:elliptical-def-intro}.
This proves \eqref{eq:effect-equivalent}.
\end{proof}

For the proof of Proposition~\ref{prop:spherical-reduction-intro}, we make use of the following spherical decomposition.
Recall that \(n_j = \dim(E_j)\) are the dimensions of the pairwise orthogonal eigenspaces \(E_1,\ldots,E_k\) associated with the Bures matrix \(K\). Further, \(\pi_j\) is the orthogonal projection onto the \(j\)th eigenspace \(E_j\) of \(K\), that is, \(K = \sum_j \kappa_j \pi_j\), where \(\kappa_j\) is the eigenvalue of \(E_j\).

\begin{lemma}[Spherical block decomposition]\label{lem:spherical-block}
Let \(U\) be uniform on \(\Sph^{d-1}\). For \(j=1,\ldots,k\), set \(D_j:=\norm{\pi_jU}^2\) and \(V_j:=\pi_jU / \norm{\pi_jU}\).
Then
\begin{enumerate}[label = (\roman*)]
    \item \label{lem:spherical-block1} \((D_1,\ldots,D_k)
\sim
\mathrm{Dirichlet}\!\left(
\frac{n_1}{2},\ldots,\frac{n_k}{2}
\right)\),
\item \label{lem:spherical-block3} \(V_j\) is uniform on the unit sphere of \(E_j\),
\item \label{lem:spherical-block2} \(V_1,\ldots,V_k\) are independent,
\item \label{lem:spherical-block4} \((V_1,\ldots,V_k)\) is independent of
\((D_1,\ldots,D_k)\).
\end{enumerate}
\end{lemma}

\begin{proof}
For each \(j\in \{1,\ldots,k\}\), choose an orthonormal basis of \(E_j\).
Since the eigenspaces are mutually orthogonal and
\[
\R^d=E_1\oplus\cdots\oplus E_k,
\]
the union of these bases forms an orthonormal basis of \(\R^d\).
Let \(G\sim N(0,I_d)\) be a standard Gaussian random vector with independent components. Then \(U := G/\norm{G}\) is uniform on \(\Sph^{d-1}\). Moreover, for \(G_j := \pi_j G\), we have
\begin{align*}
    \norm{G_j}^2\sim \chi_{n_j}^2, \qquad j=1,\ldots,k,
\end{align*}
independently.
Since \(D_j
=
\norm{\pi_jU}^2
=
\norm{G_j}^2/(\sum_{\ell=1}^k\norm{G_\ell}^2),
\)
the representation \eqref{rep_Dirichlet} with
\(\alpha_j=n_j/2\) gives
\[
(D_1,\ldots,D_k)
\sim
\mathrm{Dirichlet}\!\left(
\frac{n_1}{2},\ldots,\frac{n_k}{2}
\right).
\]
This proves \ref{lem:spherical-block1}.

For each \(j\), put \(R_j:=\norm{G_j}\) and \(V_j:=\frac{G_j}{\norm{G_j}}\),
where \(V_j\) may be defined arbitrarily on the null event
\(\{G_j=0\}\). Since \(G_j\) is a standard Gaussian vector on \(E_j\),
its distribution is invariant under orthogonal transformations of \(E_j\).
Hence \(V_j\) is uniform on the unit sphere of \(E_j\), and \(V_j\) is
independent of \(R_j\). This proves \ref{lem:spherical-block3}.

Since the random vectors \(G_1,\ldots,G_k\) are independent, also the pairs
\begin{align}\label{r_jv_jind}
    (R_1,V_1),\ldots,(R_k,V_k) \quad \text{are independent}.
\end{align}
In particular, \(V_1,\ldots,V_k\) are independent,
which proves \ref{lem:spherical-block2}.

Moreover, by \eqref{r_jv_jind} and since \(R_j\) is independent of \(V_j\) for every \(j\), we have
\[
\law\bigl((R_1,V_1),\ldots,(R_k,V_k)\bigr)
=
\bigotimes_{j=1}^k
\bigl(\law(R_j)\otimes\law(V_j)\bigr).
\]
This gives \(\law(R_1,\ldots,R_k,V_1,\ldots,V_k)
=
\law(R_1,\ldots,R_k)
\otimes
\law(V_1,\ldots,V_k)
\) and thus 
\[
(V_1,\ldots,V_k)
\perp
(R_1,\ldots,R_k).
\]
Finally, since \(D_j
=
R_j^2/ (R_1^2+\cdots+R_k^2)\) is a function of 
\((R_1,\ldots,R_k)\), we obtain
\[
(V_1,\ldots,V_k)
\perp
(D_1,\ldots,D_k),
\]
which proves \ref{lem:spherical-block4}.
\end{proof}

In the following lemma, we determine the distribution of the block-radius vector 
\[p_K(Z) = \bigl(\norm{\pi_1Z},\ldots,\norm{\pi_kZ}\bigr)\] in \eqref{def:p_K} associated with the spectral decomposition of \(K\). Recall that \(K=\sum_{j=1}^k\kappa_j\pi_j\),
where \(\pi_j\) is the orthogonal projection onto the eigenspace
\(E_j\) of \(K\). Thus \(p_K(z)
=
\bigl(\norm{\pi_1z},\ldots,\norm{\pi_kz}\bigr)
\)
records the lengths of the components of \(z\) in the distinct eigenspaces of \(K\). Its distribution for a spherical random vector depends on \(K\) only through the eigenspace dimensions \(n_j=\dim(E_j)\).

\begin{lemma}[Block-radius distribution]\label{lem:p_KZ}
    For \(Z=RU\sim\eta_\nu\), the block-radius vector in \eqref{def:p_K} satisfies
\begin{equation}\label{eq:block-radius-proof}
p_K(Z)
\eqd
\bigl(R\sqrt{D_1},\ldots,R\sqrt{D_k}\bigr),
\end{equation}
where \((D_1,\ldots,D_k)\sim \mathrm{Dirichlet}\!\left(
\frac{n_1}{2},\ldots,\frac{n_k}{2}
\right)\) is independent of \(R\).
\end{lemma}

\begin{proof}
Define \(D_j:=\norm{\pi_jU}^2\). Then, since \(Z=RU\), we have \(\norm{\pi_jZ}
=
R\norm{\pi_jU}
=
R\sqrt{D_j}\) for \(j=1,\ldots,k\).
Since \(R\) and \(U\) are independent, it follows that \((D_1,\ldots,D_k)\) is independent of \(R\), and its distribution is given by Lemma \ref{lem:spherical-block}.
\end{proof}

\begin{proof}[Proof of Proposition~\ref{prop:spherical-reduction-intro}]
By Lemma~\ref{lem:p_KZ}, we have
$\beta_i=(p_K)_\#\eta_{\nu_i}$ for $i=0,1$.
Thus, every coupling of $\eta_{\nu_0}$ and
$\eta_{\nu_1}$ induces a coupling of
$\beta_0$ and $\beta_1$.
The inequality in \eqref{eq_lowerbproj}
therefore yields
\[
W_{2,K}^2(\eta_{\nu_0},\eta_{\nu_1})
\geq
\cW_{2,\bm\kappa}^2(\beta_0,\beta_1).
\]

For the reverse inequality, let
$\gamma^*\in\Pi(\beta_0,\beta_1)$ be optimal
and consider the coupling $(Z_0^*,Z_1^*)$
constructed in \eqref{def_Zast}.
By Lemma~\ref{lem:spherical-block} and
the independence of the directions $V_j$
from the block radii, we have
\[
Z_i^*\sim\eta_{\nu_i},
\qquad i=0,1.
\]
Moreover, since both vectors share the same
unit directions within each eigenspace,
\[
(Z_0^*-Z_1^*)^\trans K(Z_0^*-Z_1^*)
=
\sum_{j=1}^k\kappa_j(r_j-s_j)^2.
\]
Consequently,
\[
\begin{aligned}
W_{2,K}^2(\eta_{\nu_0},\eta_{\nu_1})
&\leq
\E[(Z_0^*-Z_1^*)^\trans
K(Z_0^*-Z_1^*)]
=
\cW_{2,\bm\kappa}^2(\beta_0,\beta_1).
\end{aligned}
\]
Combining both inequalities completes the proof.
\end{proof}


\subsection{Proofs of Section \ref{sec_main}}

\begin{proof}[Proof of Theorem~\ref{thm:main}]
By Proposition~\ref{lem:effect-intro} and Proposition \ref{prop:spherical-reduction-intro}, we have
\[
W_2^2(P_0,P_1)
=
\norm{m_0-m_1}^2
+d_{\BW}^2(\Sigma_0,\Sigma_1)
+\cW_{2,\bm\kappa}^2(\beta_0,\beta_1),
\]
which proves \eqref{eq:main-formula}.

Now let \(\gamma^*\in\Pi(\beta_0,\beta_1)\) be optimal. Then construct
\((Z_0^*,Z_1^*)\) from \(\gamma^*\) as in \eqref{def_Zast}.
By Lemma~\ref{lem:spherical-block}, the lifted variables satisfy
\(Z_0^*\sim\eta_{\nu_0}\) and \(Z_1^*\sim\eta_{\nu_1}\).
Moreover, since \(\gamma^*\) is optimal for the reduced problem and the
lifting construction in \eqref{def_Zast} attains equality in \eqref{eq_lowerbproj},
we obtain
\begin{align}\label{eq:opt_gamma_lift}
\E[(Z_0^*-Z_1^*)^{\trans}K(Z_0^*-Z_1^*)]
&=
\cW_{2,\bm\kappa}^2(\beta_0,\beta_1)
=
W_{2,K}^2(\eta_{\nu_0},\eta_{\nu_1}),
\end{align}
where the last equality follows from \eqref{eq:spherical-spectral-reduction-intro}.

By Lemma~\ref{lem:bures-factors}, \(A_0=\Sigma_0^{1/2}\) and \(A_1=\Sigma_0^{-1/2}K\)
are covariance factors of \(\Sigma_0\) and \(\Sigma_1\), respectively.
Hence
\[
X_0=m_0+A_0Z_0^*,
\qquad
X_1=m_1+A_1Z_1^*
\]
have laws \(P_0\) and \(P_1\). Finally,
\eqref{eq:mean-separation-intro},
\eqref{eq:couplingwise-preview-intro}, and \eqref{eq:opt_gamma_lift} yield
\begin{align*}
\E\norm{X_0-X_1}^2
&=
\norm{m_0-m_1}^2
+d_{\BW}^2(\Sigma_0,\Sigma_1)
+W_{2,K}^2(\eta_{\nu_0},\eta_{\nu_1})
=
W_2^2(P_0,P_1).
\end{align*}
Thus \((X_0,X_1)\) is an optimal coupling for the quadratic cost.
\end{proof}

\begin{proof}[Proof of Corollary~\ref{cor:scalar-intro}]
Since $K=\lambda I_d$, Proposition~\ref{lem:effect-intro}
and Lemma~\ref{lem:radial-optimal} yield
\[
W_2^2(P_0,P_1)
=
\norm{m_0-m_1}^2
+d_{\BW}^2(\Sigma_0,\Sigma_1)
+\lambda W_2^2(\nu_0,\nu_1),
\]
which proves \ref{cor:scalar-intro1}.

Now assume that $\nu_0$ is non-atomic.
By Lemma~\ref{lem:radial-optimal}, the map
\[
S(z)=\frac{\tau(\norm z)}{\norm z}z,
\qquad
\tau=F_{\nu_1}^{-1}\circ F_{\nu_0},
\]
is an optimal transport from $\eta_{\nu_0}$ to
$\eta_{\nu_1}$, where the formula is understood
$\eta_{\nu_0}$-almost surely.
With $A_0=\Sigma_0^{1/2}$ and
$A_1=\Sigma_0^{-1/2}K=\lambda\Sigma_0^{-1/2}$,
transforming back to the original coordinates gives
\[
T(x)
=
m_1+A_1S\bigl(A_0^{-1}(x-m_0)\bigr),
\]
which coincides with \eqref{eq:Monge-special-intro}.
This map pushes $P_0$ to $P_1$, and its optimality
follows from \eqref{eq:couplingwise-preview-intro}
and the optimality of $S$.
This proves \ref{cor:scalar-intro2}.

Finally, define
\[
H(r):=\int_0^r\tau(s)\,\dd s,
\qquad r\geq0,
\]
using the endpoint conventions in the statement.
Since $\tau$ is nonnegative and nondecreasing,
$H$ is proper, convex, and nondecreasing.
As $t(x)=\norm{\Sigma_0^{-1/2}(x-m_0)}$ is convex,
the potential
\[
\Phi(x)=m_1^{\trans}x+\lambda H(t(x))
\]
is also proper and convex.
At points where $t(x)>0$ and $\tau$ is continuous,
the chain rule gives
\[
\nabla\Phi(x)
=
m_1+\lambda\frac{\tau(t(x))}{t(x)}
\Sigma_0^{-1}(x-m_0)
=T(x).
\]
Since $\tau$ has at most countably many
discontinuities and $\nu_0$ is non-atomic,
this identity holds $P_0$-almost surely,
proving \ref{cor:scalar-intro3}.
\end{proof}

\begin{proof}[Proof of Corollary~\ref{cor:intro-radial-bounds}]
The lower bound is a special case of Theorem \ref{thm:nested-bound}, and the upper bound is a special case of Theorem \ref{thm:multi-conditional-upper}.
\end{proof}

For the proof of Theorem \ref{thm:nested-bound}, we make use of the following lemma.

\begin{lemma}[Projected spherical laws]\label{lem:projected-deltas}
Let $Q$ be an orthogonal projection of rank $q\geq1$. For the spherical laws $\eta_i=\eta_{\nu_i}$, we have
\begin{equation}\label{eq:projection-cost-exact}
 \inf_{\pi\in\Pi(\eta_0,\eta_1)}
       \int\norm{Q(z-w)}^2\dd\pi
 =W_2^2(Q\push\eta_0,Q\push\eta_1)=\Delta_q.
\end{equation}
\end{lemma}

\begin{proof}
Put \(E:=\operatorname{ran}(Q) :=
\{Qx:x\in\R^d\}\), so that \(\dim E=q\), and define \(\mu_i^Q:=Q\push\eta_i\) for \(i=0,1\).
We first identify the projected laws. Let
\(Z_i=R_iU_i\sim\eta_i\), where \(R_i\sim\nu_i\) is independent of
\(U_i\sim\mathrm{Unif}(\Sph^{d-1})\). Applying
Lemma~\ref{lem:spherical-block} to the orthogonal decomposition \(\R^d=E\oplus E^\perp\)
shows that
\[
\|QU_i\|^2
\sim
\mathrm{Beta}\!\left(\frac q2,\frac{d-q}{2}\right)
\]
for \(q<d\), while \(\|QU_i\|=1\) if \(q=d\). Moreover, the direction
\(QU_i/\|QU_i\|\) is uniform on the unit sphere of \(E\) and independent
of \(\|QU_i\|\). Consequently, \(\mu_i^Q\) is a spherical law on \(E\)
whose radial distribution is
\[
\rho_{i,q}
:=
\law(R_i\sqrt{B_{q,d}}).
\]
Notice that \(\E\!\left[R_i^2B_{q,d}\right]
= d\,\frac qd=q,\)
so that \(\rho_{i,q}\in \cR_q\). Hence Lemma~\ref{lem:radial-optimal}, applied in
the \(q\)-dimensional space \(E\), yields
\begin{align}\label{eq:projected-radial-proof}
W_2^2(\mu_0^Q,\mu_1^Q)
=
W_2^2(\rho_{0,q},\rho_{1,q})
=
\Delta_q,
\end{align}
which proves the second equality in \eqref{eq:projection-cost-exact}.

To prove the first equality in \eqref{eq:projection-cost-exact}, let
\(\pi\in\Pi(\eta_0,\eta_1)\). Then \(\gamma:=(Q,Q)\push\pi\)
belongs to \(\Pi(\mu_0^Q,\mu_1^Q)\), and therefore
\[
\int_{\R^d\times\R^d}\|Q(z-w)\|^2\,\pi(\dd z,\dd w)
=
\int_{E\times E}\|x-y\|^2\,\gamma(\dd x,\dd y)
\geq
W_2^2(\mu_0^Q,\mu_1^Q).
\]
Taking the infimum over \(\pi\) gives
\begin{equation}\label{eq:projection-lower-proof}
\inf_{\pi\in\Pi(\eta_0,\eta_1)}
\int\|Q(z-w)\|^2\,\dd\pi
\geq
W_2^2(\mu_0^Q,\mu_1^Q).
\end{equation}

For the reverse inequality, we disintegrate \(\eta_i\) with respect to \(Q\).
Since all spaces are Polish, there exist probability kernels
\(K_i(x,\dd z)\), \(x\in E\), such that 
\[
\eta_i(\dd z)
=
\int_E K_i(x,\dd z)\,\mu_i^Q(\dd x),
\]
and \(K_i\bigl(x,\{z\in\R^d:Qz=x\}\bigr)=1\)
for \(\mu_i^Q\)-almost every \(x\). Let
\(\gamma\in\Pi(\mu_0^Q,\mu_1^Q)\) be arbitrary and define
\[
\pi_\gamma(\dd z,\dd w)
:=
\int_{E\times E}
K_0(x,\dd z)K_1(y,\dd w)\,
\gamma(\dd x,\dd y).
\]
Then \(\pi_\gamma\in\Pi(\eta_0,\eta_1)\). 
Moreover, since \(K_0\bigl(x,\{z:Qz=x\}\bigr)=1\) and \(K_1\bigl(y,\{w:Qw=y\}\bigr)=1\)
for almost every \(x\) and \(y\), Tonelli's theorem gives
\[
\begin{aligned}
\int\|Q(z-w)\|^2\,\pi_\gamma(\dd z,\dd w)
&=
\int_{E\times E}
\int
\|Qz-Qw\|^2
K_0(x,\dd z)K_1(y,\dd w)\,
\gamma(\dd x,\dd y)\\
&=
\int_{E\times E}
\|x-y\|^2\,\gamma(\dd x,\dd y).
\end{aligned}
\]
Consequently, for every
\(\gamma\in\Pi(\mu_0^Q,\mu_1^Q)\),
\begin{align*}
    \inf_{\pi\in\Pi(\eta_0,\eta_1)}
\int\|Q(z-w)\|^2\,\pi(\dd z,\dd w)
&\leq
\int\|x-y\|^2\,\gamma(\dd x,\dd y).
\end{align*}
Taking the infimum over
\(\gamma\in\Pi(\mu_0^Q,\mu_1^Q)\) therefore gives the reverse inequality. Together with \eqref{eq:projection-lower-proof}, this proves the first
equality in \eqref{eq:projection-cost-exact}. 
\end{proof}

\begin{proof}[Proof of Theorem \ref{thm:nested-bound}.]
Let $Q_j=\pi_j+\cdots+\pi_k$ be the projection onto the eigenspaces with eigenvalue at least $\kappa_j$. By \eqref{eq:spectral-layer}, the matrix can be decomposed into
\( K=  \kappa_1 I_d+\sum_{j=2}^k(\kappa_j-\kappa_{j-1})Q_j\).
For \(\pi\in\Pi(\eta_{\nu_0},\eta_{\nu_1})\), let
\((Z_0,Z_1)\sim\pi\). Then we obtain
\begin{align*}
\E_\pi\!\left[(Z_0-Z_1)^{\trans}K(Z_0-Z_1)\right]
&=
\kappa_1\E_\pi\norm{Z_0-Z_1}^2
+
\sum_{j=2}^k
(\kappa_j-\kappa_{j-1})
\E_\pi\!\left[(Z_0-Z_1)^{\trans}Q_j(Z_0-Z_1)\right] \\
&= \kappa_1\E_\pi\norm{Z_0-Z_1}^2
+
\sum_{j=2}^k
(\kappa_j-\kappa_{j-1})
\E_\pi\norm{Q_j(Z_0-Z_1)}^2,
\end{align*}
where we use for the second equality that \(
(Z_0-Z_1)^{\trans}Q_j(Z_0-Z_1)
=
\norm{Q_j(Z_0-Z_1)}^2,
\)
since each \(Q_j\) is an orthogonal projection.

Now \(I_d\) is an orthogonal projection of rank \(d\). Therefore,
Lemma~\ref{lem:projected-deltas} gives
\(\E_\pi\norm{Z_0-Z_1}^2
\geq
\Delta_d.
\)
Moreover, \(Q_j\) is the orthogonal projection onto \(E_j\oplus\cdots\oplus E_k\)
and thus has rank \(d_j=n_j+\cdots+n_k\).
Applying Lemma~\ref{lem:projected-deltas} to \(Q_j\) yields
\[
\E_\pi\norm{Q_j(Z_0-Z_1)}^2
\geq
\Delta_{d_j},
\qquad j\in \{2,\ldots,k\}.
\]
Since \(\kappa_1>0\) and \(\kappa_j-\kappa_{j-1}>0\) for \(j\), we obtain
\begin{align*}
    \E_\pi\!\left[(Z_0-Z_1)^{\trans}K(Z_0-Z_1)\right]
\geq
\kappa_1\Delta_d
+
\sum_{j=2}^k
(\kappa_j-\kappa_{j-1})\Delta_{d_j}
=
L_K.
\end{align*}
Taking the infimum over all such couplings \(\pi\) and applying \eqref{eq:Gamma-def}, gives \(\Gamma(P_0,P_1)
=
W_{2,K}^2(\eta_{\nu_0},\eta_{\nu_1}) \geq L_K\), which proves \eqref{eq:nested-bound}.

If $k\geq2$ and $\nu_0\neq\nu_1$, then
$\Delta_q>0$ for every $q=1,\ldots,d$.
Indeed, suppose that $\Delta_q=0$ for some $q\geq1$.
Then, by Lemma~\ref{lem:projected-deltas}, the orthogonal
projections of $\eta_{\nu_0}$ and $\eta_{\nu_1}$
onto any fixed $q$-dimensional subspace have
identical radial distributions.
Since these projections are spherically symmetric,
their distributions coincide. In particular, for any
unit vector $v$ in the projected subspace, we have
\(
\langle v,Z_0\rangle
\stackrel{d}{=}
\langle v,Z_1\rangle,
\)
where $Z_i\sim\eta_{\nu_i}$, $i=0,1$.
By spherical symmetry of $\eta_{\nu_0}$ and
$\eta_{\nu_1}$, all their one-dimensional projections
therefore coincide.
The Cramér--Wold theorem implies equality of the original
spherical distributions and, by taking norms,
$\nu_0=\nu_1$, a contradiction.
Consequently, since $\kappa_j-\kappa_{j-1}>0$
for $j=2,\ldots,k$, we obtain
\(L_K
=\kappa_1\Delta_d
+\sum_{j=2}^k
(\kappa_j-\kappa_{j-1})\Delta_{d_j}
>\kappa_1\Delta_d.
\)
\end{proof}

\begin{proof}[Proof of Theorem \ref{thm:multi-conditional-upper}.]
Fix $\mathbf b=(b_1,\ldots,b_m)$.
Let $e_1,\ldots,e_d$ be an orthonormal eigenbasis
of $K$, where \(e_h\) is an eigenvector associated with the eigenvalue \(\lambda_h^\downarrow\).
The subspaces
\begin{align*}
E_h:=
\operatorname{span}
\{e_{a_h+1},\ldots,e_{a_h+b_h}\},
\qquad h=1,\ldots,m.
\end{align*}
are pairwise orthogonal and have dimension $\dim E_h=b_h$. Denote by \(\pi_{E_h}= \sum_{j=a_h+1}^{a_h+b_h} e_j e_j^\top\) the orthogonal projection onto \(E_h\). Then we have
\begin{align}\label{eq:trKpiE}
\tr(K\pi_{E_h}) 
=
\tr\left(\sum_{i=1}^d
\sum_{j=a_h+1}^{a_h+b_h}
\lambda_i^\downarrow
e_i(e_i^\top e_j)e_j^\top\right)
=
\sum_{j=a_h+1}^{a_h+b_h}
\lambda_j^\downarrow 
= \Lambda_{a_h+b_h}-\Lambda_{a_h},
\end{align}
where we use for the third equality  that \(e_i\) and \(e_j\) are orthogonal.

Consider the random vectors \(Z_i = (Z_{i,1},\ldots,Z_{i,m})\) in \eqref{eq_Zi_blockradii} where \(Z_{i,h}\) models the radial part of \(\nu_i\) within the subspace \(E_h\). 
Since the Knothe--Rosenblatt coupling \((Z_0^{\mathrm{KR}},Z_1^{\mathrm{KR}})\) in \eqref{constr_KR} preserves the marginal distributions, we have 
\begin{align*}
    Z_i^{\mathrm{KR}}\eqd Z_i,
\qquad i=0,1.
\end{align*}
Let $V_1,\ldots,V_m$ be independent \(\R^d\)-valued random vectors,
independent of $(Z_0^{\mathrm{KR}},Z_1^{\mathrm{KR}})$,
where $V_h$ is uniformly distributed on the unit
sphere of $E_h$.
Define
\begin{align*}
    Y_i:=\sum_{h=1}^m Z_{i,h}^{\mathrm{KR}}V_h,
\qquad i=0,1.
\end{align*}
Note that \(Y_i\) is \(d\)-dimensional, while \(Z_{i}^{\mathrm{KR}}\) is an \(m\)-dimensional random vector.
We aim to show that
\begin{align}\label{eq_tsYKR}
    \Gamma(P_0,P_1) \leq \E\left[(Y_0-Y_1)^\top K (Y_0-Y_1)\right] = U_{K,\mathbf{b}},
\end{align}
which yields the first inequality in \eqref{thm:multi-conditional-upper1}.
Since the Knothe--Rosenblatt coupling preserves
the block-radius distributions, the spherical
block decomposition in Lemma~\ref{lem:spherical-block}
implies
\[
Y_i\sim\eta_{\nu_i},\qquad i=0,1.
\]
Thus, $(Y_0,Y_1)$ is an admissible coupling
of the standardized spherical distributions.

Now, to compute \(U_{K,\mathbf{b}}\) as an upper bound for the generator effect, we calculate
\begin{align}
    &\E\left[(Y_0-Y_1)^\top K(Y_0-Y_1)\right]
    = \sum_{h=1}^m
\E\left[
(Z_{0,h}^{\mathrm{KR}}-Z_{1,h}^{\mathrm{KR}})^2
\right]
\E[V_h^\top K V_h],
\label{eq:KR-cost-decomposition}
\end{align}
where we use \(Y_0-Y_1 = \sum_{h=1}^m (Z_{0,h}^{\mathrm{KR}}-Z_{1,h}^{\mathrm{KR}})V_h\), independence of \(Z_{i,h}^{\mathrm{KR}}\) from the directions \(V_h\), and that all mixed terms vanish because the subspaces $E_h$ are mutually orthogonal and invariant under $K$.

For the right-hand expectation in \eqref{eq:KR-cost-decomposition}, we obtain
\begin{align*}
\E[V_h^\top K V_h]
=
\tr\left(K\E[V_hV_h^\top]\right)
=
\frac{\tr(K\pi_{E_h})}{b_h}
=
\frac{\Lambda_{a_h+b_h}-\Lambda_{a_h}}{b_h},
\end{align*}
where we use for the first equality cyclicity of the trace. For the second equality, we use that $V_h$ is uniformly distributed on the
unit sphere of the $b_h$-dimensional subspace
$E_h$ and thus \(\E[V_hV_h^\top]=\frac{1}{b_h}\pi_{E_h}\). The third equality follows from \eqref{eq:trKpiE}.
It follows that
\begin{align*}
\E\left[(Y_0-Y_1)^\top K(Y_0-Y_1)\right]
=
\sum_{h=1}^m
\frac{\Lambda_{a_h+b_h}-\Lambda_{a_h}}{b_h}
\E\left[
(Z_{0,h}^{\mathrm{KR}}-Z_{1,h}^{\mathrm{KR}})^2
\right]
=
\sum_{h=1}^m
\frac{\Lambda_{a_h+b_h}-\Lambda_{a_h}}{b_h}
C_h^{\mathbf b}
=U_{K,\mathbf b}.
\end{align*}
Since $(Y_0,Y_1)$ is an admissible coupling
of $\eta_{\nu_0}$ and $\eta_{\nu_1}$, we have proven \eqref{eq_tsYKR}.

Finally, choosing $m=1$ and $\mathbf b=(d)$,
the Knothe--Rosenblatt coupling reduces to the
comonotone coupling of the radial laws \(\nu_0\) and \(\nu_1\). 
Hence, \(C_1^{(d)} = W_2^2(\nu_0,\nu_1)=\Delta_d = \delta_R^2\).
Since, by definition in \eqref{sum_eigenvalue}, $\Lambda_d=\tr K$, it follows that 
\[U_K \leq U_{K,(d)} = \frac{\tr K}{d}\delta_R^2,\]
which proves the second inequality in
\eqref{thm:multi-conditional-upper1}
and completes the proof.
\end{proof}

\subsection{Proofs of Section \ref{sec:dyn}}

For the proof of Theorem \ref{lem:KR-student-one-dimensional}, we determine the conditional distributions \(Q_{a,b}\) in  Table \ref{tab:KR-families} of the block radius vectors \(Z_i = (Z_{i,1},\ldots,Z_{i,m})\) in \eqref{eq_Zi_blockradii}.

\begin{lemma}[Conditional distributions of block radii]
\label{lem:KR-conditional-block-radii}
For $\mathbf b=(b_1,\ldots,b_m)$ with integers \(b_h\geq 1\) such that \(\sum_h b_h = d\), let $Z=(Z_{1},\ldots,Z_{m})$ be a block-radius vector as in \eqref{eq_Zi_blockradii}. 
Assume that $Z$ is the block-radius vector of
a spherical random vector belonging to one
of the four families in Table~\ref{tab:KR-families}, with associated constants $c>0$ and 
$\varepsilon\in\{-1,0,1\}$. Set $a_h=\sum_{\ell<h}b_\ell$.
Then, for $h=1,\ldots,m$, the conditional
distribution of $Z_{h}$ given the preceding
block radii $Z_{1:h-1}=r$ satisfies
\begin{align}\label{eq:KR-conditional-radius}
\left(Z_h/\alpha_{h}\mid Z_{1:h-1}=r\right)
\sim \, Q_{a_h,b_h}, \qquad \alpha_{h} := \sqrt{c+\varepsilon\sum_{\ell<h}r_\ell^2}
\end{align}
where the distribution $Q_{a_h,b_h}$ specified in Table \ref{tab:KR-families} does not depend on \(r\).
\end{lemma}

\begin{proof}
We distinguish the four distributions considered in Table \ref{tab:KR-families}.

\emph{Gaussian distributions.}
Let $G\sim N(0,I_d)$ and decompose
$G=(G_1,\ldots,G_m)$ according to the
ordered composition $\mathbf b=(b_1,\ldots,b_m)$,
where $G_h\in\R^{b_h}$.
Since the covariance matrix is $I_d$, the
random vectors $G_1,\ldots,G_m$ are independent
with $G_h\sim N(0,I_{b_h})$.
Consequently, the block radii satisfy
\[
Z_h=\|G_h\|=\sqrt{H_h},
\qquad H_h:=\|G_h\|^2\sim\chi^2_{b_h},
\]
where $H_1,\ldots,H_m$ are independent.
In particular, for almost every
$r=(r_1,\ldots,r_{h-1})$, we obtain
\begin{align*}
(Z_h\mid Z_{1:h-1}=r)
&\eqd Z_h
=\sqrt{H_h}
\sim Q_{a_h,b_h}.
\end{align*}
Since $c=1$ and $\varepsilon=0$, this gives \eqref{eq:KR-conditional-radius}.

\emph{Student-$t$ distributions.}
Recall that
a covariance-standardized spherical Student-$t$
distribution with $\theta>2$ degrees of freedom
admits the representation
$\sqrt{c/V}\,G$, where $c=\theta-2$,
$G\sim N(0,I_d)$ and $V\sim\chi^2_\theta$
are independent.
Decomposing $G$ into independent Gaussian blocks
$G_h\sim N(0,I_{b_h})$, $h=1,\ldots,m$,
the block-radius vector $Z=(Z_1,\ldots,Z_m)$
in \eqref{eq_Zi_blockradii} satisfies
\begin{align}\label{eq:KR-student-block-rep}
Z\eqd
\sqrt{\frac{c}{V}}\,
\left(\|G_1\|,\ldots,\|G_m\|\right)
\eqd
\left(
\sqrt{\frac{cH_1}{V}},
\ldots,
\sqrt{\frac{cH_m}{V}}
\right),
\end{align}
where $H_h:=\|G_h\|^2\sim\chi^2_{b_h}$.
Then, conditional on $V=v$, the block radii
$Z_1,\ldots,Z_m$ are independent.
In particular, the joint conditional density
of the first $h-1$ block radii, evaluated at
$r=(r_1,\ldots,r_{h-1})$, is given by
\begin{align*}
&f_{Z_{1:h-1}\mid V=v}(r)
=\prod_{\ell<h}f_{Z_\ell\mid V=v}(r_\ell)
\propto
v^{a_h/2}
\exp\left(
-\frac{v}{2 c}\sum_{\ell<h}r_\ell^2
\right),
\end{align*}
where $a_h:=\sum_{\ell<h}b_\ell$.
Here, $\propto$ means proportional as a function
of $v$ (i.e., $v$-independent factors are omitted).
Multiplying the above likelihood by the density of
$V\sim\chi^2_{\theta}$ and applying Bayes'
formula yields
\begin{align*}
f_{V\mid Z_{1:h-1}=r}(v)
&\propto
v^{a_h/2}
\exp\left(
-\frac{v}{2c}\sum_{\ell<h}r_\ell^2
\right)
\frac{1}{2^{\theta/2}\Gamma(\theta/2)}
v^{\theta/2-1}e^{-v/2}\\
&\propto
v^{(\theta+a_h)/2-1}
\exp\left[
-\frac v2\left(
1+\frac{\sum_{\ell<h}r_\ell^2}{c}
\right)
\right].
\end{align*}
This gives 
\begin{align*}
(V\mid Z_{1:h-1}=r)
\eqd
\frac{c}{c+\sum_{\ell<h}r_\ell^2}
\chi^2_{\theta+a_h}.
\end{align*}
Since, by \eqref{eq:KR-student-block-rep},
$Z_h=\sqrt{c H_h/V}$,
where $H_h\sim\chi^2_{b_h}$ is independent
of $V$ and the preceding block radii
$Z_{1:h-1}$, we obtain
\begin{align*}
(Z_h\mid Z_{1:h-1}=r)
\eqd
\sqrt{c+ \sum_{\ell<h}r_\ell^2}\,
\sqrt{\frac{H_h}{V_{\theta+a_h}}},
\end{align*}
where $V_{\theta+a_h}\sim\chi^2_{\theta+a_h}$
is independent of $H_h$. Choosing \(\varepsilon = 1\), this establishes \eqref{eq:KR-conditional-radius}.

\emph{Uniform distributions on balls.}
Suppose that \(Y=(Y_1,\ldots,Y_d)\) has the spherical law uniform on
the ball of radius $\sqrt{d+2}$. Let $c=d+2$.
Let $G\sim N(0,I_d)$ and $E\sim\chi^2_2$
be independent. We first observe that
\begin{align}\label{eq:uniform-ball-gaussian-rep}
Y\eqd
\sqrt{\frac{c}{\|G\|^2+E}}\,G.
\end{align}
Indeed, since $\|G\|^2\sim\chi^2_d$, we have \(\|G\|^2/(\|G\|^2+E)
\sim\mathrm{Beta}\left(d/2,1\right)\)
Consequently, the radius of the random vector
on the right-hand side of
\eqref{eq:uniform-ball-gaussian-rep} satisfies
\begin{align}\label{eq_unitsphe}
\Pp\left(
\sqrt{\frac{c\|G\|^2}{\|G\|^2+E}}\leq s
\right)
=
\left(\frac{s}{\sqrt c}\right)^d,
\qquad 0\leq s\leq\sqrt c.
\end{align}
Moreover, its direction $G/\|G\|$ is uniform
on $\mathbb S^{d-1}$ and independent of its
radius. This establishes
\eqref{eq:uniform-ball-gaussian-rep}.

Now decompose $G=(G_1,\ldots,G_m)$
according to the ordered composition
$\mathbf b=(b_1,\ldots,b_m)$, where
$G_h\sim N(0,I_{b_h})$ are independent.
Setting
\[
H_h:=\|G_h\|^2\sim\chi^2_{b_h},
\qquad h=1,\ldots,m,
\]
we obtain from \eqref{eq:uniform-ball-gaussian-rep}
that the block-radius vector \(Z=(Z_1,\ldots,Z_m)\), defined as in \eqref{eq_Zi_blockradii}, admits the representation
\begin{align*}
Z
&\eqd
\sqrt{\frac{c}{H_1+\cdots+H_m+E}}
\left(\sqrt{H_1},\ldots,\sqrt{H_m}\right).
\end{align*}
In particular, 
\begin{align}\label{eq_rep_zh}
\frac{Z_h^2}{c}
\eqd
\frac{H_h}{H_1+\cdots+H_m+E},
~ h=1,\ldots,m,
\qquad \text{and} \qquad
1-\frac{\sum_{\ell=1}^m Z_\ell^2}{c}
\eqd
\frac{E}{H_1+\cdots+H_m+E}.
\end{align}
Since $H_1,\ldots,H_m,E$ are independent
chi-squared random variables with degrees
of freedom $b_1,\ldots,b_m,2$, respectively,
we have
\begin{align*}
&\left(
\frac{Z_1^2}{c},\ldots,\frac{Z_m^2}{c},
1-\frac{\sum_{\ell=1}^m Z_\ell^2}{c}
\right)\sim
\mathrm{Dirichlet}
\left(
\frac{b_1}{2},\ldots,\frac{b_m}{2},1
\right),
\end{align*}
see \eqref{rep_Dirichlet}.
To determine the conditional distribution of \(Z_h\mid Z_{1:h-1}=r\) let
\begin{align*}
W_h:=\sum_{\ell=h}^m H_\ell+E,
\qquad h=1,\ldots,m, \qquad W_{m+1}:= E.
\end{align*}
Using the representation \eqref{eq_rep_zh}, we obtain
\begin{align*}
Z_h^2\eqd c\frac{H_h}{W_1},
\qquad
c-\sum_{\ell<h}Z_\ell^2
\eqd c\frac{W_h}{W_1}, \qquad \text{and thus} \quad 
\frac{Z_h^2}
{c-\sum_{\ell<h}Z_\ell^2}
\eqd
\frac{H_h}{W_h}
=
\frac{H_h}{H_h+W_{h+1}}.
\end{align*}
Since $H_h\sim\chi^2_{b_h}$ and \(W_{h+1}\sim\chi^2_{d-a_h-b_h+2}\)
are independent, the relationship between the beta and the chi-squared distribution in \eqref{eq_betachi} yields
\[
\frac{H_h}{H_h+W_{h+1}}
\sim
\mathrm{Beta}\left(
\frac{b_h}{2},
\frac{d-a_h-b_h}{2}+1\right) \qquad \text{and} \qquad \frac{H_h}{H_h+W_{h+1}} \perp W_h,
\]
where we use that \(W_h = H_h+W_{h+1}\) is chi-squared distributed as a sum of independent chi-squared distributed random variables.
Since \(H_1,\ldots,H_m,E\) are independent, the random vector \((H_1,\ldots,H_{h-1})\) is independent of \((H_h,W_{h+1})\) and, consequently, independent of \(H_h/W_h\).
Therefore, for almost every $r=(r_1,\ldots,r_{h-1})$, we obtain
\begin{align*}
&\left(
\frac{Z_h^2}
{c-\sum_{\ell<h}Z_\ell^2}
\,\middle|\,Z_{1:h-1}=r
\right)\sim
\mathrm{Beta}\left(
\frac{b_h}{2},
\frac{d-a_h-b_h}{2}+1
\right)
\end{align*}
independent of \(r\).
Choosing $\varepsilon=-1$, this establishes \eqref{eq:KR-conditional-radius}
for uniform distributions on balls.

\emph{Uniform distributions on spheres.}
Suppose that $Y=(Y_1,\ldots,Y_d)$ has the spherical
law uniform on the sphere of radius $\sqrt d$, and let $c=d$.
Let $G\sim N(0,I_d)$. Since $G/\|G\|$ is
uniformly distributed on the unit sphere $\mathbb S^{d-1}$, setting \(E=0\),
we have as in \eqref{eq:uniform-ball-gaussian-rep} that
\begin{align}\label{eq:uniform-sphere-gaussian-rep}
Y\eqd\sqrt{\frac{d}{\|G\|^2 + E}}\,G \qquad \text{and} \qquad \left\| \sqrt{\frac{c}{\left\|G\right\|^2 +E }}\, G\right\|^2 = \frac{c \,\lVert G \rVert^2}{\lVert G\rVert^2 + E} = c.
\end{align}
Now, proceeding as in the case of the uniform distribution on the ball, we obtain \(W_{h+1} \sim \chi^2_{d-a_h-b_h}\) whenever \(a_h+b_h< d\) (instead of \(W_{h+1}\sim \chi^2_{d-a_h-b_h+2}\)). Then we obtain
\begin{align*}
    \left(
\frac{Z_h^2}
{c-\sum_{\ell<h}Z_\ell^2}
\,\middle|\,Z_{1:h-1}=r
\right)
\eqd
\frac{H_h}{H_h+W_{h+1}}
\sim
\mathrm{Beta}\left(
\frac{b_h}{2},
\frac{d-a_h-b_h}{2}
\right)
\end{align*}
for almost every $r=(r_1,\ldots,r_{h-1})$. If \(a_h+b_h=d\), then \(W_{h+1}=0\) and the above conditional distribution is the Dirac distribution in \(1\).
Choosing \(\varepsilon = -1\), this establishes \eqref{eq:KR-conditional-radius} for uniform distributions on spheres.
\end{proof}

\begin{proof}[Proof of Theorem \ref{lem:KR-student-one-dimensional}.]
Let $\mathbf b=(b_1,\ldots,b_m)$ be an ordered composition of $d$. We first show that the Knothe--Rosenblatt coupling in \eqref{constr_KR} admits for the spherical distributions in Table \ref{tab:KR-families} the cross-moment identity representation
\begin{align}\label{eq:student-KR-crossmoment}
\E\left[
Z_{0,h}^{\mathrm{KR}}Z_{1,h}^{\mathrm{KR}}
\right]
=
\sqrt{c_0c_1}\,
I_{a_h,b_h}
\prod_{\ell<h}J_{a_\ell,b_\ell},
\qquad h=1,\ldots,m.
\end{align}
To prove this identity, we use Lemma \ref{lem:KR-conditional-block-radii} to obtain for the random vector of block radii \(Z_i = (Z_{i,1},\ldots,Z_{i,m})\) in \eqref{eq_Zi_blockradii} the conditional distributions
\begin{align*}
    \left(\frac{Z_{i,h}}{\sqrt{c_i+\varepsilon_i\sum_{\ell<h}r_\ell^2}}\mid Z_{i,1:h-1}=r\right)
\sim Q_{i;a_h,b_h},
\end{align*}
Thus, the conditional quantile construction
of the Knothe--Rosenblatt coupling in
\eqref{constr_KR} yields
\begin{align}\label{eq:student-KR-recursion}
Z_{i,h}^{\mathrm{KR}}
=
\sqrt{c_i+
\varepsilon_i \sum_{\ell<h}(Z_{i,\ell}^{\mathrm{KR}})^2}
\,g_{i;a_h,b_h}(U_h).
\end{align}
This implies 
\begin{align*}
c_i+\varepsilon_i\sum_{\ell\leq h}(Z_{i,\ell}^{\mathrm{KR}})^2
=
\left(
c_i+\varepsilon_i\sum_{\ell<h}(Z_{i,\ell}^{\mathrm{KR}})^2
\right)
\left(1+\varepsilon_i \,g_{i;a_h,b_h}(U_h)^2\right).
\end{align*}
Iterating the above identity yields
\begin{align}\label{eq:student-KR-recursion2}
\sqrt{c_i+\varepsilon_i\sum_{\ell<h}(Z_{i,\ell}^{\mathrm{KR}})^2}
=
\sqrt{c_i}
\prod_{\ell<h}
\sqrt{1+\varepsilon_i\, g_{i;a_\ell,b_\ell}(U_\ell)^2}.
\end{align}
Hence, combining \eqref{eq:student-KR-recursion} and \eqref{eq:student-KR-recursion2} gives
\begin{align}\label{rep_produ_KR}
Z_{i,h}^{\mathrm{KR}}
=
\sqrt{c_i}\,
g_{i;a_h,b_h}(U_h)
\prod_{\ell<h}
\sqrt{1+\varepsilon_i\, g_{i;a_\ell,b_\ell}(U_\ell)^2}.
\end{align}
Since $U_1,\ldots,U_m$ are independent, we obtain
\begin{align}\nonumber
\E\left[
Z_{0,h}^{\mathrm{KR}}Z_{1,h}^{\mathrm{KR}}
\right]
&=
\sqrt{c_0c_1}\,
\E\left[
g_{0;a_h,b_h}(U_h)g_{1;a_h,b_h}(U_h)
\right]
\prod_{\ell<h}
\E\left[
\sqrt{1+\varepsilon_0\, g_{0;a_\ell,b_\ell}(U_\ell)^2}
\sqrt{1+\varepsilon_1\, g_{1;a_\ell,b_\ell}(U_\ell)^2}
\right]
\\
&=
\sqrt{c_0c_1}\,
I_{a_h,b_h}
\prod_{\ell<h}J_{a_\ell,b_\ell},\label{eq_misch}
\end{align}
which proves \eqref{eq:student-KR-crossmoment}.

Next, we show the representation \eqref{eq:student-KR-cost} below.
Since \(Z_i^{\mathrm{KR}} \eqd Z_i\) and \(Z_{i,h} = R_i\sqrt {D_h}\) for \(D_h\sim \mathrm{Beta}(b_h/2,(d-b_h)/2)\) independent of \(R_i\), we obtain \(\E[(Z_{i,h}^{\mathrm{KR}})^2] = \E[Z_{i,h}^2] = \E[R_i^2] \E [D_h] = d \tfrac {b_h}d = b_h\). Hence, we obtain with \eqref{eq_misch} that
the quadratic cost of the $h$th block is
\begin{align*}
C_h^{\mathbf b}
&=
\E\left[
(Z_{0,h}^{\mathrm{KR}}
-Z_{1,h}^{\mathrm{KR}})^2
\right]
=
2b_h
-2\sqrt{c_0c_1}\,
I_{a_h,b_h}
\prod_{\ell<h}J_{a_\ell,b_\ell}.
\end{align*}
Using \(w_{a,b}:=(\Lambda_{a+b}-\Lambda_a)/b\), we obtain
\[
\sum_{h=1}^m b_hw_{a_h,b_h}
=
\sum_{h=1}^m
(\Lambda_{a_h+b_h}-\Lambda_{a_h})
=
\tr K.
\]
Hence, by definition of $U_{K,\mathbf b}$ in \eqref{def_UkB}, we have
\begin{align}\label{eq:student-KR-cost}
U_{K,\mathbf b}
= \sum_{h=1}^m w_{a_h,b_h}C_h^{\mathbf b}=
2\tr K
-2\sqrt{c_0c_1}
\sum_{h=1}^m
w_{a_h,b_h}I_{a_h,b_h}
\prod_{\ell<h}J_{a_\ell,b_\ell}.
\end{align}

It remains to minimize
\eqref{eq:student-KR-cost} over all ordered
compositions $\mathbf b = (b_1,\ldots,b_m)$ of $d$.
Since $c_0,c_1>0$, this is equivalent to maximizing
\[
S(\mathbf b)
:=
\sum_{h=1}^m
w_{a_h,b_h}I_{a_h,b_h}
\prod_{\ell<h}J_{a_\ell,b_\ell}.
\]
We show by backward induction that the recursion \eqref{eq:student-KR-DP} yields the desired maximum.
For \(a=d\), we set \(B_d = 0\) and there is nothing to optimize. For the induction step \(a \mapsto a-1\), assume that the assertion holds for all \(B_j\) with \(j\geq a\).
Consider an arbitrary ordered composition
$\mathbf b=(b_1,\ldots,b_m)$ of $d-a+1$, and set
\[
a_h:=a-1+\sum_{\ell<h}b_\ell.
\]
Separating the contribution of the block with index \(h=1\)
and factoring out $J_{a-1,b_1}$, we obtain
\begin{align*}
\sum_{h=1}^m
w_{a_h,b_h}I_{a_h,b_h}
\prod_{\ell<h}J_{a_\ell,b_\ell}
=w_{a-1,b_1}I_{a-1,b_1}
+
J_{a-1,b_1}
\sum_{h=2}^m
w_{a_h,b_h}I_{a_h,b_h}
\prod_{\ell=2}^{h-1}J_{a_\ell,b_\ell},
\end{align*}
using the convention that the empty product \(\prod_{\ell < 1} J_{a_\ell,b_\ell}\) equals \(1\).
By the induction hypothesis, the maximum
of the right-hand sum over all ordered compositions $(b_2,\ldots,b_m)$ is
$B_{a-1+b_1}$ because \(a-1+b_1\geq a\).
Hence, we maximize
\[
w_{a-1,b_1}I_{a-1,b_1}
+J_{a-1,b_1}B_{a-1+b_1}
\]
over all possible first block sizes
$1\leq b_1\leq d-a+1$. This yields
\begin{align*}
    B_{a-1} = \max_{1\leq b\leq d-a+1} \left\{ w_{a-1,b}I_{a-1,b} +J_{a-1,b}B_{a-1 +b} \right\}.
\end{align*}
This establishes the induction step. 
Consequently, $B_0$ is the maximum of
$S(\mathbf b)$ over all ordered compositions
of $d$. Therefore,
\[
U_K
=
\min_{\mathbf b}U_{K,\mathbf b}
=
2\tr K-2\sqrt{c_0c_1}\,B_0,
\]
which completes the proof.
\end{proof}

\subsection{Proofs of Section \ref{sec:simulation}}

\begin{proof}[Proof of Lemma \ref{lem:normalized-projection-costs}.]
Let \(E\subset\R^d\) be an \(r\)-dimensional subspace, and choose an
orthonormal basis \(e_1,\ldots,e_r\) of \(E\). For every subset
\(A\subset\{1,\ldots,r\}\) with \(|A|=q\), let \(Q_A\) denote the
orthogonal projection onto \(E_A:=\operatorname{span}\{e_a:a\in A\}\).
Let $\mu_i^E$ be the projection of
$\eta_{\nu_i}$ onto $E$, $i=0,1$, and let
\(\gamma\in\Pi(\mu_0^E,\mu_1^E)\). By
\eqref{eq:projection-cost-exact},
\begin{align*}
    \int\|Q_A(x-y)\|^2\,\gamma(\dd x,\dd y) \geq \Delta_q
\end{align*}
for every such \(A\). Summing over all
\(\binom rq\) subsets \(A\) gives 
\begin{align}\label{est_gamma}
\binom rq\Delta_q
&\leq
\sum_{\substack{A\subset\{1,\ldots,r\}\\|A|=q}}
\int\|Q_A(x-y)\|^2\,\gamma(\dd x,\dd y)
=
\binom{r-1}{q-1}
\int\|x-y\|^2\,\gamma(\dd x,\dd y),
\end{align}
where we use for the second equality that 
\begin{align*}
    \sum_{\substack{A\subset\{1,\ldots,r\}\\|A|=q}} \lVert Q_A v\rVert^2 
    = \sum_{\substack{A\subset\{1,\ldots,r\}\\|A|=q}}\sum_{\ell\in A} \langle v,e_\ell\rangle^2 
    = \binom{r-1}{q-1} \sum_{\ell = 1}^r \langle v,e_\ell\rangle^2
    = \binom{r-1}{q-1} \lVert v \rVert^2.
\end{align*}
Hence, \eqref{est_gamma} gives
\[
\int\|x-y\|^2\,\dd\gamma
\geq
\frac{\binom rq}{\binom{r-1}{q-1}}\Delta_q
=
\frac rq\,\Delta_q.
\]
Taking the infimum over \(\gamma\) yields \(\Delta_r\geq\frac rq\Delta_q\) and thus \(\Delta_q/q\leq \Delta_r/r\).
\end{proof}

\section*{Acknowledgements}

The first author acknowledges support from the Austrian
Science Fund (FWF) [10.55776/PAT1669224] through the
project \emph{Stochastic orders for functional dependence}.
During the preparation of this manuscript, the authors
used ChatGPT (OpenAI) for assistance with mathematical
arguments, technical checks, literature searches,
numerical computations, and language editing.
The authors take full responsibility for the
mathematical content and conclusions of this work.

\bibliographystyle{plainnat}
\bibliography{elliptical_wasserstein_arXiv}

@article{BhatiaJainLim2019,
  author  = {Bhatia, R. and Jain, T. and Lim, Y.},
  title   = {On the {Bures--Wasserstein} distance between positive definite matrices},
  journal = {Expositiones Mathematicae},
  volume  = {37},
  number  = {2},
  pages   = {165--191},
  year    = {2019},
  doi     = {10.1016/j.exmath.2018.01.002}
}

@article{Brenier1987,
  author  = {Brenier, Y.},
  title   = {D{\'e}composition polaire et r{\'e}arrangement monotone des champs de vecteurs},
  journal = {Comptes Rendus de l'Acad{\'e}mie des Sciences de Paris, S{\'e}rie I},
  volume  = {305},
  number  = {19},
  pages   = {805--808},
  year    = {1987}
}

@article{Brenier1991,
  author  = {Brenier, Y.},
  title   = {Polar factorization and monotone rearrangement of vector-valued functions},
  journal = {Communications on Pure and Applied Mathematics},
  volume  = {44},
  number  = {4},
  pages   = {375--417},
  year    = {1991},
  doi     = {10.1002/cpa.3160440402}
}

@article{CambanisHuangSimons1981,
  author  = {Cambanis, S. and Huang, S. and Simons, G.},
  title   = {On the theory of elliptically contoured distributions},
  journal = {Journal of Multivariate Analysis},
  volume  = {11},
  number  = {3},
  pages   = {368--385},
  year    = {1981},
  doi     = {10.1016/0047-259X(81)90082-8}
}

@article{ChenZhang2024,
  author        = {Chen, K. and Zhang, Y.},
  title         = {Optimal transport and {Wasserstein} barycenter for radially contoured distributions},
  journal       = {Analysis Mathematica},
  volume        = {51},
  number        = {3},
  pages         = {775--792},
  year          = {2025},
  doi           = {10.1007/s10476-025-00105-5},
  eprint        = {2404.08383},
  archivePrefix = {arXiv}
}

@article{CuestaRuschendorfTuero1993,
  author  = {Cuesta-Albertos, J. A. and R{\"u}schendorf, L. and Tuero-D{\'i}az, A.},
  title   = {Optimal coupling of multivariate distributions and stochastic processes},
  journal = {Journal of Multivariate Analysis},
  volume  = {46},
  number  = {2},
  pages   = {335--361},
  year    = {1993},
  doi     = {10.1006/jmva.1993.1064}
}

@article{DowsonLandau1982,
  author  = {Dowson, D. C. and Landau, B. V.},
  title   = {The {Fr\'echet} distance between multivariate normal distributions},
  journal = {Journal of Multivariate Analysis},
  volume  = {12},
  number  = {3},
  pages   = {450--455},
  year    = {1982},
  doi     = {10.1016/0047-259X(82)90077-X}
}

@book{FangKotzNg1990,
  author    = {Fang, K.-T. and Kotz, S. and Ng, K. W.},
  title     = {Symmetric Multivariate and Related Distributions},
  publisher = {Chapman and Hall},
  address   = {London},
  series    = {Monographs on Statistics and Applied Probability},
  volume    = {36},
  year      = {1990}
}

@article{Gelbrich1990,
  author  = {Gelbrich, M.},
  title   = {On a formula for the {$L^2$} {Wasserstein} metric between measures on {Euclidean} and {Hilbert} spaces},
  journal = {Mathematische Nachrichten},
  volume  = {147},
  number  = {1},
  pages   = {185--203},
  year    = {1990},
  doi     = {10.1002/mana.19901470121}
}

@article{GhaffariWalker2021,
  author  = {Ghaffari, N. and Walker, S. G.},
  title   = {{Parseval}'s identity and optimal transport maps},
  journal = {Statistics \& Probability Letters},
  volume  = {170},
  pages   = {108989},
  year    = {2021},
  doi     = {10.1016/j.spl.2020.108989}
}

@article{GhaffariWalker2023,
  author  = {Ghaffari, N. and Walker, S. G.},
  title   = {{$W_2$} barycenters for radially related distributions},
  journal = {Statistics \& Probability Letters},
  volume  = {195},
  pages   = {109788},
  year    = {2023},
  doi     = {10.1016/j.spl.2023.109788}
}

@article{GivensShortt1984,
  author  = {Givens, C. R. and Shortt, R. M.},
  title   = {A class of {Wasserstein} metrics for probability distributions},
  journal = {Michigan Mathematical Journal},
  volume  = {31},
  number  = {2},
  pages   = {231--240},
  year    = {1984},
  doi     = {10.1307/MMJ/1029003026}
}

@article{GunasingamWong2025,
  author  = {Gunasingam, M. and Wong, T.-K. L.},
  title   = {Adapted optimal transport between {Gaussian} processes in discrete time},
  journal = {Electronic Communications in Probability},
  volume  = {30},
  pages   = {1--14},
  year    = {2025},
  doi     = {10.1214/25-ECP654}
}

@article{HammKorzeniowski2024,
  author  = {Hamm, K. and Korzeniowski, A.},
  title   = {On {Wasserstein} distances for affine transformations of random vectors},
  journal = {Foundations of Data Science},
  volume  = {6},
  number  = {4},
  pages   = {468--491},
  year    = {2024},
  doi     = {10.3934/fods.2024023}
}

@book{KotzBalakrishnanJohnson2000,
  author    = {Kotz, S. and Balakrishnan, N. and Johnson, N. L.},
  title     = {Continuous Multivariate Distributions, Volume 1: Models and Applications},
  edition   = {2},
  publisher = {Wiley},
  year      = {2000},
  doi       = {10.1002/0471722065}
}

@book{MarshallOlkinArnold2011,
  author    = {Marshall, A. W. and Olkin, I. and Arnold, B. C.},
  title     = {Inequalities: Theory of Majorization and Its Applications},
  edition   = {2},
  series    = {Springer Series in Statistics},
  publisher = {Springer},
  address   = {New York},
  year      = {2011},
  doi       = {10.1007/978-0-387-68276-1}
}

@inproceedings{MuzellecCuturi2018,
  author        = {Muzellec, B. and Cuturi, M.},
  title         = {Generalizing point embeddings using the {Wasserstein} space of elliptical distributions},
  booktitle     = {Advances in Neural Information Processing Systems},
  volume        = {31},
  pages         = {10258--10269},
  year          = {2018},
  eprint        = {1805.07594},
  archivePrefix = {arXiv}
}

@article{OlkinPukelsheim1982,
  author  = {Olkin, I. and Pukelsheim, F.},
  title   = {The distance between two random vectors with given dispersion matrices},
  journal = {Linear Algebra and its Applications},
  volume  = {48},
  pages   = {257--263},
  year    = {1982},
  doi     = {10.1016/0024-3795(82)90112-4}
}

@book{RachevRuschendorf1998,
  author    = {Rachev, S. T. and R{\"u}schendorf, L.},
  title     = {Mass Transportation Problems, Volume {I}: Theory},
  series    = {Probability and Its Applications},
  publisher = {Springer},
  address   = {New York},
  year      = {1998}
}

@article{Ruschendorf1985,
  author  = {R{\"u}schendorf, L.},
  title   = {The {Wasserstein} distance and approximation theorems},
  journal = {Zeitschrift f{\"u}r Wahrscheinlichkeitstheorie und Verwandte Gebiete},
  volume  = {70},
  pages   = {117--129},
  year    = {1985},
  doi     = {10.1007/BF00532240}
}

@article{RuschendorfRachev1990,
  author  = {R{\"u}schendorf, L. and Rachev, S. T.},
  title   = {A characterization of random variables with minimum {$L^2$}-distance},
  journal = {Journal of Multivariate Analysis},
  volume  = {32},
  number  = {1},
  pages   = {48--54},
  year    = {1990},
  doi     = {10.1016/0047-259X(90)90070-X}
}

@misc{SeiMatsuzoe2025,
  author        = {Sei, T. and Matsuzoe, H.},
  title         = {Open problems in information geometry: a discussion at {FDIG} 2025},
  howpublished  = {arXiv preprint arXiv:2509.06989},
  year          = {2025},
  eprint        = {2509.06989},
  archivePrefix = {arXiv},
  url           = {https://arxiv.org/abs/2509.06989}
}

@article{Takatsu2011,
  author  = {Takatsu, A.},
  title   = {{Wasserstein} geometry of {Gaussian} measures},
  journal = {Osaka Journal of Mathematics},
  volume  = {48},
  number  = {4},
  pages   = {1005--1026},
  year    = {2011},
  doi     = {10.18910/4973}
}

@book{Villani2009,
  author    = {Villani, C.},
  title     = {Optimal Transport: Old and New},
  publisher = {Springer},
  address   = {Berlin},
  series    = {Grundlehren der mathematischen Wissenschaften},
  volume    = {338},
  year      = {2009},
  doi       = {10.1007/978-3-540-71050-9}
}

@article{FournierGuillin2015,
  author  = {Fournier, N. and Guillin, A.},
  title   = {On the rate of convergence in Wasserstein
             distance of the empirical measure},
  journal = {Probability Theory and Related Fields},
  volume  = {162},
  pages   = {707--738},
  year    = {2015},
  doi     = {10.1007/s00440-014-0583-7}
}

@article{WeedBach2019,
  author  = {Weed, J. and Bach, F.},
  title   = {Sharp asymptotic and finite-sample rates
             of convergence of empirical measures
             in Wasserstein distance},
  journal = {Bernoulli},
  volume  = {25},
  number  = {4A},
  pages   = {2620--2648},
  year    = {2019},
  doi     = {10.3150/18-BEJ1065}
}

@inproceedings{Cuturi2013,
  author    = {Cuturi, M.},
  title     = {Sinkhorn Distances: Lightspeed
               Computation of Optimal Transport},
  booktitle = {Advances in Neural Information
               Processing Systems},
  volume    = {26},
  year      = {2013}
}

@inproceedings{GenevayEtAl2019,
  author    = {Genevay, A. and Chizat, L.
               and Bach, F. and Cuturi, M.
               and Peyr{\'e}, G.},
  title     = {Sample Complexity of Sinkhorn Divergences},
  booktitle = {Proceedings of the 22nd International
               Conference on Artificial Intelligence
               and Statistics},
  series    = {Proceedings of Machine Learning Research},
  volume    = {89},
  pages     = {1574--1583},
  year      = {2019}
}

@incollection{KuhnEtAl2019,
  author    = {Kuhn, D. and Mohajerin Esfahani, P. and Nguyen, V. A. and Shafieezadeh-Abadeh, S.},
  title     = {{Wasserstein} distributionally robust optimization: Theory and applications in machine learning},
  booktitle = {Operations Research \& Management Science in the Age of Analytics},
  series    = {INFORMS TutORials in Operations Research},
  publisher = {INFORMS},
  pages     = {130--166},
  year      = {2019},
  doi       = {10.1287/educ.2019.0198}
}

@article{AmariMatsuda2024,
  author  = {Amari, S.-I. and Matsuda, T.},
  title   = {Information geometry of {Wasserstein} statistics on shapes and affine deformations},
  journal = {Information Geometry},
  volume  = {7},
  number  = {2},
  pages   = {285--309},
  year    = {2024},
  doi     = {10.1007/s41884-024-00139-y}
}

@article{ChenWangZhang2026,
  author  = {Chen, K. and Wang, Z. and Zhang, Y.},
  title   = {A Relaxed {Wasserstein} Distance Formulation
             for Mixtures of Radially Contoured Distributions},
  journal = {Communications on Applied Mathematics
             and Computation},
  year    = {2026},
  doi     = {10.1007/s42967-026-00579-6}
}

\end{document}